\documentclass[12pt]{amsart}
\usepackage{amsmath,amssymb,bm,color,mathrsfs}
\usepackage{graphicx}
\usepackage{slashbox,multirow}
\usepackage{algorithm}
\usepackage{float}
\begin{document}
\newtheorem{theorem}{Theorem}[section]
\newtheorem{lemma}[theorem]{Lemma}
\newtheorem{definition}[theorem]{Definition}
\newtheorem{corollary}[theorem]{Corollary}
\newtheorem{example}[theorem]{Example}
\newtheorem{remark}[theorem]{Remark}
\numberwithin{equation}{section}
\numberwithin{table}{section}
\numberwithin{figure}{section}
\numberwithin{algorithm}{section}
\allowdisplaybreaks[4]
\def\red{\color{red}}
\def\p{\partial}
\def\O{\Omega}
\def\R{\mathbb{R}}
\def\LT{{L^2(\O)}}
\def\pli{\psi_\ell^{\omega_i}}
\def\plpi{\psi_{\ell'}^{\omega_i}}
\def\Pil{\Phi_{i,\ell}}
\def\til{\phi_{i,\ell}}
\def\cT{\mathcal{T}}
\def\cN{\mathcal{N}}
\def\PiR{\Pi_{\scriptscriptstyle R}}
\def\zVh{\mathaccent'27{V}_h}
\def\fC{\mathfrak{C}}
\def\cCk{\fC_{h,k}}
\def\KER{\mathrm{Ker}\,\Pi_{aux}}
\def\MS{{V}_{aux}^{\raise 2pt\hbox{$\scriptstyle\rm {ms},h$}}}
\def\MSk{{V}_{aux,k}^{\raise 2pt\hbox{$\scriptstyle\rm {ms},h$}}}
\def\MSu{{u}_{aux}^{\mathrm{ms},h}}
\def\MSku{{u}_{aux,k}^{\mathrm{ms},h}}
\def\MSL{V_{aux,k}^{\hspace{1pt}\lower 3pt\hbox{$\scriptstyle\rm {ms},h$}}}
\def\MSLu{u_{aux,k}^{\mathrm{ms},h}}
\def\d{\displaystyle}
\def\sC{\mathfrak{C}}
\def\patch{{\omega_j}}
\def\cV{\mathcal{V}}
\def\tV{\tilde{V}}
\def\tv{\tilde{v}}
\def\cI{\mathcal{I}}
\def\tV{\tilde{V}}
\def\sS{\mathscr{S}}
\def\bA{\mathbf{A}}
\def\bv{\mathbf{v}}
\def\bw{\mathbf{w}}
\def\bx{\mathbf{x}}
\def\bK{\mathbf{K}}
\def\bC{\mathbf{C}}
\def\bB{\mathbf{B}}
\def\bI{\mathbf{I}}
\def\hp{\hat{p}}
\def\hphi{\hat{\phi}}
\def\Vd{V_h^{d}}
\def\mystrut(#1,#2){\vrule height #1pt depth #2pt width 0pt}
\def\defeq{\,\raise 7pt\hbox{\tiny def}\hspace{-13pt}=\,}
%
\title{A Spectral LOD Method for Multiscale Problems with High Contrast}
\thanks{This  work  was supported in part by
the National Science Foundation under
 Grant No. DMS-22-08404 and Grant No. DMS-25-13273.}
\author{Susanne C. Brenner}
\address{Susanne C. Brenner, Department of Mathematics and Center for
Computation and Technology, Louisiana State University, Baton Rouge,
LA 70803, USA}
\email{brenner@math.lsu.edu}
\author{Jos\'e C. Garay}
\address{Jos\'e C. Garay, Institute of Mathematics, Augsburg University,
Augsburg, Germany}\email{jocgafer@gmail.com}
\author{Li-yeng Sung}
\address{Li-yeng Sung,
 Department of Mathematics and Center for Computation and Technology,
 Louisiana State University, Baton Rouge, LA 70803, USA}
\email{sung@math.lsu.edu}
\begin{abstract}
  We present a multiscale finite element method for a diffusion problem
  with rough and  high contrast coefficients.  The construction of the
  multiscale finite element space
  is based on the localized orthogonal decomposition methodology and it involves
  solutions of local finite element eigenvalue problems.  We show that the
  performance of
  the multiscale finite element method is similar to the performance of
  standard finite element methods
  {for the  homogeneous Dirichlet boundary value problem for the
  Poisson equation on smooth or convex domains.}  Simple explicit error
  estimates are established
  under  conditions that can be verified from the outputs of the computation.
\end{abstract}
\date{November 1, 2025}
\maketitle
\section{Introduction}\label{sec:Introduction}
 Let $\O$ be the unit square $(0,1)\times(0,1)$.
 We consider the following model problem: Find $u\in H^1_0(\O)$ such that
\begin{equation}\label{eq:BVP}
 a(u,v)=\int_\Omega fv\,dx\qquad \forall\,v\in H^1_0(\O),
\end{equation}
  where
\begin{equation*}
  a(u,v)=\int_\O \kappa\nabla u\cdot\nabla v\,dx,
\end{equation*}
 and the diffusion coefficient
 $\kappa(x)\in L^\infty(\O)$ is a piecewise constant function such that $\min_{x\in\O} \kappa(x)=1$ and
 $\max_{x\in\O}\kappa(x)=\beta$.
\begin{remark}\label{rem:Notation}\rm
We follow the standard notation for function spaces and norms that
can be found for example in \cite{ADAMS:2003:Sobolev}.
\end{remark}
\par
 Accordingly, we have
\begin{equation}\label{eq:Comparison}
  |v|_{H^1(\O)}\leq \|v\|_a \qquad\forall\,v\in H^1(\O),
\end{equation}
 where $\|v\|_a=\sqrt{a(v,v)}$ is the energy norm.
\par
\begin{figure}[H]
\begin{center}
 \includegraphics[height=2in]{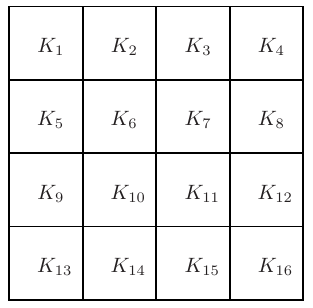} \qquad
 \includegraphics[height=2in]{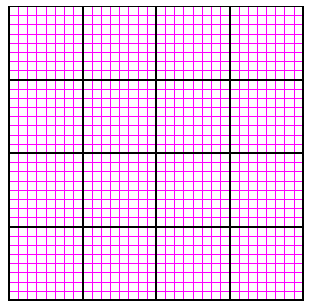}
\end{center}
 \caption{A coarse mesh $\cT_H$ with 16 elements (left) and a fine mesh $\cT_h$ (right).}
 \label{fig:Meshes}
 \end{figure}

  Let $\cT_H$ be a uniform coarse triangulation of $\O$ with ${m}$ (closed) elements
  $K_1,\ldots,K_m$
  (cf. Figure~\ref{fig:Meshes} (left)), where $H$ is
  the length of the edges of the squares in $\cT_H$. The $Q_1$ finite element subspace of
  $H^1_0(\O)$ associated with $\cT_H$ (cf. \cite{Ciarlet:1978:FEM,BScott:2008:FEM})
   is denoted by $V_H$.
  \par
  Let
  $\cT_h$ be a fine triangulation of $\O$ obtained by a uniform refinement of $\cT_H$
  (cf. Figure~\ref{fig:Meshes} (right)).
  The $Q_1$ finite element subspace of $H^1_0(\O)$ associated with $\cT_h$
  is denoted by
  $V_h$.
  We use $\cV_h$ to denote the set of {interior} vertices of $\cT_h$.
\par\smallskip
  Let $u_h\in V_h$ be the standard finite element solution of \eqref{eq:BVP}, i.e.,
\begin{equation}\label{eq:uh}
  a(u_h,v)=\int_\O fv\,dx \qquad\forall\,v\in V_h.
\end{equation}
\begin{remark}\label{rem:uhBdd}\rm
  We have a Poincar\'e-Friedrichs inequality
\begin{equation}\label{eq:PF}
  \|v\|_\LT\leq |v|_{H^1(\O)} \qquad\forall\,v\in H^1_0(\O)
\end{equation}
 on the unit square $\O$, which together with \eqref{eq:Comparison} and \eqref{eq:uh}
 implies that
\begin{equation}\label{eq:uhBdd}
  \|u_h\|_a\leq \|f\|_\LT.
\end{equation}
\end{remark}
  We assume $u_h$
  is a good approximation of $u$, which means in particular that $\kappa(x)$ should be
  resolved by $\cT_h$.
    However
  the computation of $u_h$ is expensive because $\text{dim}\,V_h\gg1$, especially if we have
  to solve for different $f$'s repeatedly.
  On the other hand the standard finite element solution $u_H\in V_H$, while
  affordable because $\text{dim}\,V_H\ll \text{dim}\,V_h$, is a poor approximation of
 $u$.  Therefore we need to design a multiscale finite element method that bridges the two scales.
\par
{
 We will construct the multiscale finite element method through the localized
 orthogonal decomposition
 strategy (cf. \cite{HP:2013:LOD,MP:2014:LOD,MP:2021:LOD})
  and by using a coarse space consisting  of local eigenfunctions.
\par
  There is a substantial  literature on solving multiscale problems by using local eigenfunctions
  (cf. \cite{GE:2010:DDMS1,GE:2010:DDMS2,BL:2011:Multiscale,SDHNPS:2014:EigenDD,KC:2015:Spectral,
  CEL:2018:Multiscale,MS:2021:Spectral,MS:2022:Error,MSD:2022:GFEM,MS:2025:Spectral}).
  Our approach is related to the work in  \cite{CEL:2018:Multiscale}. The global multiscale
  finite element
  spaces in \cite{CEL:2018:Multiscale} and this paper are identical, albeit their bases
  are generated differently.
  On the other hand, the localized multiscale finite element spaces are different.
  Instead of solving
  a local version of the global problem as in \cite{CEL:2018:Multiscale}, we construct
  our multiscale
  finite element space $\MSk$ by solving the corrector equation through $k$ iterations of
  a preconditioned conjugate gradient algorithm.
  Our main theoretical result is an estimate of the form
\begin{equation*}
  \|u_h-\MSLu\|_a\leq C_\dag H \|f\|_\LT,
\end{equation*}
 provided that $k\geq k_\dag$.  Here $\MSLu$ is the Galerkin solution of \eqref{eq:BVP}
 for the space $\MSk$ and the explicit constants $C_\dag$ and $k_\dag$ are calculated from
 the outputs of the computation in the construction of $\MSk$.
  {Thus the construction and analysis of our method for \eqref{eq:BVP} is  {\em a posteriori}
  in nature.  In particular
    the analysis does not require the
 exponential decay of the basis functions of the multiscale finite element space.}
}
\par
{
 The rest of the paper is organized as follows.  We introduce the auxiliary
 space $V_{aux}$ comprised of
 local eigenfunctions in Section~\ref{sec:Auxiliary}, and construct in Section~\ref{sec:PiAux}
 an operator $\Pi_{aux}$ from $V_h$ to $V_{aux}$ that has a desired approximation property.
 The ideal (or global) multiscale finite element method is analyzed in Section~\ref{sec:Ideal},
 followed
 by the construction of a localized multiscale finite element method in Section~\ref{sec:LOD}.
 The analysis of the localized multiscale method is carried out in Section~\ref{sec:Analysis} and
 corroborating numerical results are presented in Section~\ref{sec:Numerics}.  We end the paper
 with some concluding remarks in Section~\ref{sec:Conclusions}.} {Appendix~\ref{append:EVs}
 provides some eigenvalue estimates that are used in the computation.}
%
\section{The Auxiliary Space $V_{aux}$}\label{sec:Auxiliary}
 On each element $K_i \in\cT_H$, we construct the local auxiliary space $V_{aux}^{(i)}$
 by solving an eigenvalue problem.
\par
 Let $V_h(K_i)$ be the restriction of $V_h$ to $K_i$ and $n_i=\mathrm{dim}V_h(K_i)$.
 Note that $v\in V_h(K_i)$
 vanishes on $\p K_i\cap\p\O$.   We will
 denote by $\zVh(K_i)$ the subspace of
 $V_h(K_i)$ whose members vanish on $\p K_i$.
\par
 According to the finite dimensional Spectral Theorem (cf. \cite[Section~79]{Halmos:1942:FDVS}),
 there exist eigenvalues $0\leq\lambda_1^{(i)}\leq \cdots\leq \lambda_{n_i}^{(i)}$
 and corresponding eigenfunctions $\psi_1^{(i)},\ldots,\psi_{n_i}^{(i)}$ such that
\begin{align}
  \int_{K_i}\kappa\nabla \psi_j^{(i)}\cdot\nabla w\,dx&=
  \lambda_j^{(i)} H^{-2}\int_{K_i} \kappa
  \psi_j^{(i)} \,w\,dx  \qquad\forall\,w\in V_h(K_i),\label{eq:EigenEquation}\\
  H^{-2}\int_{K_i}\kappa (\psi_j^{(i)})^2dx&=1.\label{eq:Normalization}
\end{align}
 \par
 We can rewrite \eqref{eq:EigenEquation} and \eqref{eq:Normalization} concisely as
\begin{align*}
 a_i(\psi_j^{(i)},w)= \lambda_j^{(i)}s_i(\psi_j^{(i)},w)\quad\forall\,w\in V_h(K_i)
 \quad\text{and}\quad
 s_i(\psi_j^{(i)},\psi_j^{(i)})=1,
\end{align*}
 where
\begin{alignat*}{3}
 a_i(v,w)&=\int_{K_i}\kappa \nabla v\cdot\nabla w\,dx&\qquad&\forall\,v,w\in V_h(K_i),\\
 s_i(v,w)&=H^{-2}\int_{K_i}\kappa vw\,dx&\qquad&\forall\,v,w\in V_h(K_i).
\end{alignat*}
\par
 The functions $\psi_j^{(i)}$ for $1\leq j\leq n_i$ form an orthonormal basis of the
 inner product space $\big(V_h(K_i),s_i(\cdot,\cdot)\big)$.
 Given any $v\in V_h(K_i)$, we can write
\begin{equation}\label{eq:Representation}
  v=\sum_{j=1}^{n_i} c_j\psi_j^{(i)}
\end{equation}
 and then we have
\begin{align}
  \|v\|_{L^2(K_i;\kappa)}^2&=\int_{K_i}\kappa v^2dx=H^2\sum_{j=1}^{n_i}c_j^2,
  \label{eq:LocalLTwo}\\
  \|v\|_{a_i}^2&=\int_{K_i} \kappa|\nabla v|^2dx=\sum_{j=1}^{n_i}\lambda_j^{(i)}c_i^2.
  \label{eq:LocalEnergyNorm}
\end{align}
\par
 We define
\begin{equation}\label{eq:LocalAuxSpace}
 V_{aux}^{(i)}=\text{span}\{\psi_1^{(i)},\ldots,\psi_{L_i}^{(i)}\},
\end{equation}
 where $L_i$ is sufficiently large so that
\begin{equation}\label{eq:Choice}
  \lambda_{j} ^{(i)}>\frac12 \mu_{i} \quad\text{if}\quad j>L_i,
\end{equation}
 and $\mu_{i}$ is the first nonzero eigenvalue of the analog
  of \eqref{eq:EigenEquation}--\eqref{eq:Normalization}
 for the Laplace operator (i.e., $\kappa=1$) on $K_i$.
{
\begin{remark}\label{rem:Scaling}\rm
  Due to the scaling factor $H^{-2}$ in  \eqref{eq:EigenEquation},
  the eigenvalue $\mu_i$ in \eqref{eq:Choice}
  is $\approx 1$.  Moreover it can be shown (cf. Appendix~\ref{append:EVs}) that
\begin{alignat*}{3}
 \mu_i&>\pi^2&\qquad&\text{if $K_i$ is disjoint from $\p\O$},\\
 \mu_i&>\pi^2/4&\qquad&\text{if only one of the edges of $K_i$ is on $\p\O$},\\
 \mu_i&>\pi^2/2&\qquad&\text{if two of the edges of $K_i$ are on $\p\O$}.
\end{alignat*}
\end{remark}
}
\begin{remark}\label{rem:VAuxDimension}\rm
 $L_i$ is bounded
 by the number of high contrast channels or inclusions in $K_i$
 (cf. \cite[Appendix~A]{GE:2010:DDMS2}).
\end{remark}
\begin{remark}\label{rem:LocalEigenCost}\rm
  The cost for computing $\psi_1^{(i)},\ldots,\psi_{L_i}^{(i)}$ is $O((n/m)^3)$,
  where $n=\mathrm{dim}\,V_h$,
  and the eigenfunctions associated with
   different elements of $\cT_H$ can be computed in parallel.
\end{remark}
\par
 The global auxiliary space
\begin{equation*}
 V_{aux}=V_{aux}^{(1)}\oplus\cdots\oplus V_{aux}^{(m)}
\end{equation*}
  is the direct sum of the local auxiliary spaces, and
\begin{equation}\label{eq:DimensionVaux}
  \mathrm{dim}V_{aux}=\sum_{i=1}^m \mathrm{dim}V_{aux}^{(i)}=\sum_{i=1}^m L_i \defeq L.
\end{equation}
%
%
\section{The operator $\Pi_{aux}$}\label{sec:PiAux}
 Let $\d v=\sum_{j=1}^{n_i}c_j\psi_j^{(i)}$ be an arbitrary function in $V_h(K_i)$.
 The local projection operator $\Pi_{aux}^{(i)}:V_h(K_i)\longrightarrow V_{aux}^{(i)}$
 is defined by
\begin{equation}\label{eq:LocalPiAux}
 \Pi_{aux}^{(i)}v=\sum_{j=1}^{L_i}c_j\psi_j^{(i)}
 =\sum_{j=1}^{L_i} s_i(v,\psi_j^{(i)})\psi_j^{(i)},
\end{equation}
 i.e., $\Pi_{aux}^{(i)}$ is the orthogonal projection operator associated with the inner product
 $s_i(\cdot,\cdot)$.
\par
 It follows from \eqref{eq:Representation}, \eqref{eq:LocalEnergyNorm} and
 \eqref{eq:LocalPiAux} that
\begin{equation}\label{eq:LocalPiAuxEnergyEstimate}
  \|\Pi_{aux}^{(i)}v\|_{a_i}^2=\sum_{1\leq j\leq L_i}\lambda_j^{(i)}c_i^2
  \leq \sum_{1\leq j\leq n_i} \lambda_j^{(i)}c_i^2=\|v\|_{a_i}^2,
\end{equation}
 and in view of \eqref{eq:LocalLTwo}, \eqref{eq:LocalEnergyNorm}, \eqref{eq:Choice}
 and Remark~\ref{rem:Scaling},
\begin{align}\label{eq:LocalPiAuxLTwoEstimate}
  H^{-2}\|v-\Pi_{aux}^{(i)}v\|_{L^2(K_i;\kappa)}^2
  &=\sum_{L_i<j\leq n_i}c_j^2\\
  &\leq \sum_{L_i<j\leq n_i}[\lambda_j^{(i)}]^{-1}\lambda_j^{(i)}c_j^2
  \leq \frac{2}{\mu_i}\sum_{L_i<j\leq n_i}\lambda_j^{(i)}c_j^2
  \leq C_*^2\|v\|_{a_i}^2,\notag
\end{align}
{
 where (cf. Remark~\ref{rem:Scaling})
\begin{equation}\label{eq:C*}
 C_*=2^{3/2}\pi \geq \big(2\max_{1\leq i\leq m}(1/\mu_i)\big)^\frac12.
\end{equation}
}
\par
 Let 
\begin{equation*}
\tilde V_h=V_h(K_1)\oplus\cdots\oplus V_h(K_m)
\end{equation*}
be the direct sum of $V_h(K_i)$.
 The projection operator $\Pi_{aux}:\tilde V_h\longrightarrow V_{aux} \,(\subset \tV_h)$
 is given by
\begin{equation*}
 \big(\Pi_{aux}\tv\big)_i=\Pi_{aux}^{(i)}v_i,
\end{equation*}
 where $\tv=(v_1,\ldots,v_m)\in\tV_h$ (with $v_i\in V_h(K_i)$).
\par
 For a function $\tv\in \tV_h$, we define
\begin{align}
 \|\tv\|_{L^2(\O;\kappa)}&=\Big(\sum_{i=1}^m \int_{K_i}\kappa v_i^2dx\Big)^\frac12
 =\Big(\sum_{i=1}^m \|v_i\|_{L^2(K_i;\kappa)}^2\Big)^\frac12,\label{eq:LTNorm}\\
 \|\tv\|_a&=\Big(\sum_{i=1}^m \int_{K_i}\kappa|\nabla v_i|^2dx\Big)^\frac12
 =\Big(\sum_{i=1}^m \|v_i\|_{a_i}^2\Big)^\frac12.\label{eq:EnergyNorm}
\end{align}
\par
 The  following result is an immediate consequence of the local estimates
 \eqref{eq:LocalPiAuxEnergyEstimate}, \eqref{eq:LocalPiAuxLTwoEstimate}
 and the definitions \eqref{eq:LTNorm} and \eqref{eq:EnergyNorm}.
\begin{lemma}\label{lem:PiAuxApproximation}
  We have
\begin{equation*}
  \|\Pi_{aux}\tv\|_a\leq \|\tv\|_a  \qquad \forall\,\tv\in\tV_h
\end{equation*}
 and
\begin{equation*}
 \|\tv-\Pi_{aux}\tv\|_{L^2(\O;\kappa)}\leq {C_*}H\|\tv\|_a
 \qquad \forall\,\tv\in\tV_h.
\end{equation*}
\end{lemma}
\begin{remark}\label{rem:Vh}\rm
  We can identify $v\in V_h$ with $\tv=(v_1,\ldots,v_m)\in\tV_h$, where
  $v_i=v\big|_{K_i}$ is the restriction
  of $v$ on $K_i$,
  and we have
\begin{align*}
 \|\tv\|_{L^2(\O;\kappa)}=\Big(\int_\O \kappa v^2dx\Big)^\frac12
 =\|v\|_{L^2(\O;\kappa)} \quad\text{and}\quad
 \|\tv\|_a=\Big(\int_\O \kappa |\nabla v|^2dx\Big)^\frac12=\|v \|_a.
\end{align*}
\end{remark}
\par
 According to Remark~\ref{rem:Vh}, we can identify $V_h$ with a subspace of $\tV_h$. Hence
 $\Pi_{aux}$ is defined on $V_h$, and
\begin{equation}\label{eq:PiAuxVh}
 \Pi_{aux}v=\big(\Pi_{aux}^{(i)}\big(v\big|_{K_1}\big),\ldots,\Pi_{aux}^{(m)}
 \big(v\big|_{K_m}\big)\big) \qquad\forall\,v\in V_h.
\end{equation}
 But the restriction of $\Pi_{aux}$ to $V_h$, still denoted by $\Pi_{aux}$,
 is {no longer} a projection
 because $V_{aux}\subset \tV_h$ is not a subspace of $V_h\subset\tV_h$.
{
\begin{remark}\label{rem:Maier}\rm
  The idea of using a nonconforming coarse subspace in the construction of a
  localized orthogonal decomposition
  method was also employed in \cite{Maier:2021:HigherOrder}.
\end{remark}
}
\par
 We will use $\KER$ to represent the kernel of $\Pi_{aux}$ in $V_h$, i.e.,
\begin{equation}\label{eq:Kernel}
 \KER=\{v\in V_h:\,\Pi_{aux}v=0\}=\{v\in V_h:\,\Pi_{aux}^{(i)} \big(v\big|_{K_i}\big)=0
  \quad\text{for}\quad 1\leq i\leq m\}.
\end{equation}
\par
 The functions in $\KER$ satisfy an important  estimate that follows from
 Lemma~\ref{lem:PiAuxApproximation}
 and Remark~\ref{rem:Vh}:
 \begin{equation}\label{eq:KernelEst}
   \|w\|_{L^2(\O;\kappa)}=\|w-\Pi_{aux}w\|_{L^2(\O;\kappa)}\leq {C_*} H \|w\|_a
   \qquad\forall\,w\in\KER.
 \end{equation}
%
\section{The ideal (or global) multiscale finite element method}\label{sec:Ideal}
 The ideal (or global) multiscale finite element space $\MS$ is the subspace of
  $V_h$ orthogonal
 to $\KER$ with respect to the bilinear form $a(\cdot,\cdot)$, i.e.,
\begin{equation}\label{eq:IdealMS}
 \MS=\{v\in V_h:\,a(v,w)=0\quad \forall\,w\in\KER\}.
\end{equation}
 It follows that
\begin{equation}\label{eq:DimensionMS}
  \mathrm{dim}\MS=\mathrm{dim}V_h-\mathrm{dim}\,\KER.
\end{equation}
\par
 The ideal multiscale finite element method is to find $\MSu\in\MS$ defined by
\begin{equation}\label{eq:IdealMsFEM}
 a(\MSu,v)=\int_\O fv\,dx\qquad\forall\,v\in\MS.
\end{equation}
\par
 There is a simple error analysis for the ideal multiscale finite element method.
\begin{theorem}\label{thm:IdealErrors}
 Let $u_h\in V_h$ be the solution of \eqref{eq:uh}.
 We have
\begin{equation}\label{eq:IdealEnergyError}
 \|{u_h}-\MSu\|_a\leq {C_*}H\|\kappa^{-1/2}f\|_\LT\leq {C_*}H\|f\|_\LT,
\end{equation}
where ${C_*}$ is the positive constant defined by \eqref{eq:C*}.
\end{theorem}
\begin{proof} From \eqref{eq:IdealMS} and the Galerkin orthogonality
 $$a(u_h-\MSu,w)=0 \qquad\forall\,w\in \MS$$
that follows from \eqref{eq:BVP} and \eqref{eq:IdealMsFEM}, we see that $u_h-\MSu\in\KER$.
\par
 Therefore we have, by \eqref{eq:BVP}, \eqref{eq:KernelEst}, \eqref{eq:IdealMS}
  and the Cauchy-Schwarz inequality,
\begin{align*}
 \|u_h-\MSu\|_a^2&=a(u_h-\MSu,u_h-\MSu)\\
   &=a(u_h,u-\MSu)\\
      &=\int_\O f(u_h-\MSu)dx\\
         &\leq \|\kappa^{-1/2}f\|_\LT\|u_h-\MSu\|_{L^2(\O:\kappa)}\\
          &\leq \|\kappa^{-1/2}f\|_\LT C_*H\|u_h-\MSu\|_a,
\end{align*}
 which implies \eqref{eq:IdealEnergyError}.
\end{proof}
\par
 It follows from Theorem~\ref{thm:IdealErrors} and a standard duality argument that
\begin{equation}\label{eq:IdealLTError}
 \|u_h-\MSu\|_{L^2(\O;\kappa)}\leq (C_*H)^2\|\kappa^{-1/2}f\|_\LT\leq (C_*H)^2\|f\|_\LT.
\end{equation}
\subsection{A Basis for $\KER$}\label{subsec:Ideal}
 In order to construct a basis for $\MS$,
 first we need a basis for $\KER$.
\par
\begin{figure}[htbp]
\includegraphics[width=1in]{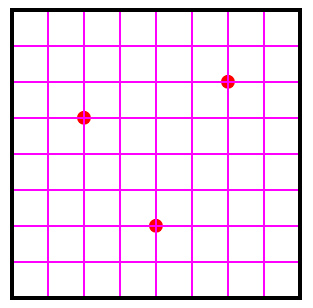}
\caption{{The dual nodes $\hp_1^{(i)}$, $\hp_2^{(i)}$ and
$\hp_3^{(i)}$ in $K_i$ represented by
the red} {dots, where $L_i=3$.}}
\label{fig:DualNodes}
\end{figure}

\smallskip\noindent
{\bf Assumption} \quad Given $K_i$, we can
find $L_i$ nodes $\hp_1^{(i)},\ldots,\hp_{L_i}^{(i)}$ of
 $\cT_h$ in the interior of $K_i$ (cf. Figure~\ref{fig:DualNodes} where $L_i=3$)
 such that the matrix $S_i$ defined by
\begin{equation*}
 S_i(j,k)=s_i(\hphi_j^{(i)},\psi_k^{(i)}) \qquad 1\leq j,k\leq L_i
\end{equation*}
 is nonsingular, where $\hphi_j^{(i)}$ is the hat function (from $V_h$) associated
 with the node $\hp_j$ normalized so that
\begin{equation}\label{eq:NormalizedHatFunction}
 a_i(\hphi_j^{(i)},\hphi_j^{(i)})=1.
\end{equation}
 We will refer to these nodes as the dual nodes and denote by $\Vd$ the
 space spanned by the functions
 $\hphi_j^{(i)}$ for $1\leq i \leq m$ and $1\leq j\leq L_i$.   In view of
 \eqref{eq:DimensionVaux}, we have
\begin{equation}\label{eq:dimVhd}
  \mathrm{dim} V_h^d=L=\mathrm{dim} V_{aux}.
\end{equation}
\begin{remark}\label{rem:AssumptionII}\rm
 As $h$ goes to zero, the eigenfunctions $\psi_j^{(i)}$ ($1\leq j\leq L_i$)
  approach the eigenfunctions
 of the corresponding continuous eigenproblem posed on the subspace of $H^1(K_i)$
  whose members vanish on $\p K_i\cap\p\O$,
 and hence they can be approximated accurately in the norm of $L^2(K_i;\kappa)$ by
 functions in $\zVh(K_i)$.
 Consequently, for $h$ sufficiently small,
  there exist functions $\zeta_1,\ldots,\zeta_{L_i}\in\zVh(K_i)$ such that
 the matrix $\sS_i$ defined by
\begin{equation*}
  \sS_i(j,k)=s_i(\zeta_j,\psi_k^{(i)}) \qquad 1\leq j,k\leq L_i
\end{equation*}
 is nonsingular.  This means the rank of the matrix with components $s_i(\phi_\ell^{(i)},\psi_k^{(i)})$
 for $1\leq \ell\leq \tilde n_i$ and $1\leq k\leq L_i$ is $L_i$,  where $\tilde n_i=\mathrm{dim}\zVh(K_i)$ and
 $\phi_\ell^{(i)}$ ($1\leq\ell\leq\tilde n_i$) is the normalized hat
  function associated with the node
 $p_\ell$ interior to $K_i$.  Therefore the existence of dual nodes is guaranteed.
   Moreover in this case
   it is also very likely that a set of $L_i$ randomly
  chosen interior nodes of $K_i$ can  be used as dual nodes.
\end{remark}
\begin{remark}\label{rem:ChoosingNodes}\rm
 In practice we just choose the dual nodes $\hp_1,\ldots,\hp_{L_i}$ randomly and verify that
  the  matrix $S_i$ is nonsingular.  It is also convenient for the error analysis that
  the dual nodes are chosen to be separated, i.e., any element of $\cT_h$ can have
  at most one of the dual nodes
  as a vertex (cf. Figure~\ref{fig:DualNodes}).  Consequently we have
\begin{equation}\label{eq:DualKronecker}
  a(\hphi_j^{(i)},\hphi_{j'}^{(i')})=\begin{cases}
    1 &\qquad\text{if $i=i'$ and $j=j'$}\\[2pt]
    0&\qquad\text{otherwise}
  \end{cases}.
\end{equation}
\end{remark}
\par
 Observe that the existence of dual nodes in $K_i$ implies that
 $\Pi_{aux}^{(i)}:\zVh(K_i)\longrightarrow V_{aux}^{(i)}$
 is a surjection.  It follows that $\Pi_{aux}:V_h\longrightarrow V_{aux}$
 is also a surjection
 and hence (cf. \eqref{eq:DimensionVaux})
\begin{equation}\label{eq:KernelDimension}
  \mathrm{dim}\,\KER=\mathrm{dim}V_h-\mathrm{dim}V_{aux}=\mathrm{dim}V_h-L,
\end{equation}
 which together with \eqref{eq:DimensionMS} implies that
\begin{equation}\label{eq:ExplicitDimensionMS}
  \mathrm{dim}\MS=L.
\end{equation}
\par
 Let the functions $\tilde\phi_1^{(i)},\ldots,\tilde\phi_{L_i}^{(i)}$ in
 $\mathrm{span}\{\hphi_1^{(i)},\ldots,\hphi_{L_i}^{(i)}\}$ satisfy
\begin{equation}\label{eq:Kronecker}
  s_i(\tilde\phi_j^{(i)},\psi_k^{(i)})=\delta_{jk}=\begin{cases}
  1&\quad\text{if $j=k$}\\[4pt]
  0&\quad\text{if $j\neq k$}
  \end{cases}.
\end{equation}
 These dual functions  are given by
\begin{equation*}
 \tilde\phi_j^{(i)}=\sum_{\ell=1}^{L_i} \tau_{j\ell}^{(i)}\hphi_\ell^{(i)} \qquad
 \text{for}\quad 1\leq j\leq L_i,
\end{equation*}
 where $\tau_{j\ell}^{(i)}$ are the components of the matrix $S_i^{-1}$.
 Note that
\begin{equation}\label{eq:SameSpan}
 \mathrm{span}\{\tilde\phi_1^{(i)},\ldots,\tilde\phi_{L_i}^{(i)}\}=
 \mathrm{span}\{\hphi_1^{(i)},\ldots,\hphi_{L_i}^{(i)}\}.
\end{equation}
\par
{
 We define the number
\begin{equation}\label{eq:Mi}
  M_i=\text{the 2-norm of the $L_i\times L_i$ matrix with components
  $a(\tilde\phi_j^{(i)},\tilde\phi_\ell^{(i)})$}.
\end{equation}
 It will play a role in the error analysis in Section~\ref{sec:Analysis}.}
\begin{remark}\label{rem:DualFunctions}\rm
 It follows from \eqref{eq:LocalPiAux} and \eqref{eq:Kronecker} that
 \begin{equation}\label{eq:DualFunctionProjecion}
 \Pi_{aux}^{(i)}\tilde\phi_j^{(i)}=\psi_j^{(i)} \qquad\text{for}\quad 1\leq j\leq L_i.
 \end{equation}
  Consequently $\Pi_{aux}:V_h^d\longrightarrow V_{aux}$ is an
  isomorphism by \eqref{eq:DimensionVaux},
  \eqref{eq:dimVhd} and \eqref{eq:SameSpan}.
\end{remark}
\par
\begin{figure}[htbp]
\includegraphics[height=1.6in]{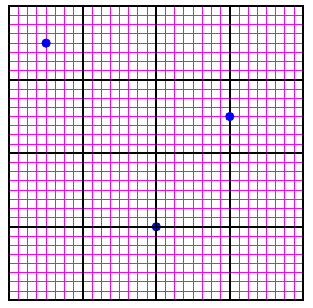}
\caption{The three possibilities for $p\in \cV_h\setminus\cN_{aux}$}
\label{fig:Locations}
\end{figure}

  Let $\cN_{aux}=\{\hp_j^{(i)}:\,1\leq i\leq m, 1\leq j\leq L_i\}$ be the set of the dual nodes
  in $\cT_h$.
  Given any node $p\in \cV_h\setminus\cN_{aux}$, denote by $\phi_{p}$ the hat function
  (in $V_h$) associated with $p$ normalized by
 $\|\phi_p\|_a=1$.
  Let the function $\tilde\phi_p\in V_h$ be defined by
\begin{equation}\label{eq:tildephip}
  \tilde\phi_p=\sum_{i=1}^m \sum_{j=1}^{L_i} s_i(\phi_p,\psi_j^{(i)})\tilde\phi_j^{(i)}.
\end{equation}
\begin{lemma}\label{lem:KernalBasis}
 The function $\phi_p-\tilde\phi_p$ belongs to $\KER$ and
 $\{\phi_p-\tilde\phi_p:\,p\in\cV_h\setminus\cN_{aux}\}$ is a basis of $\KER$.
\end{lemma}
\begin{proof} It follows from \eqref{eq:LocalPiAux}, \eqref{eq:DualFunctionProjecion}
and  \eqref{eq:tildephip} that
\begin{equation*}
 \Pi_{aux}^{(i)}\phi_p=\sum_{j=1}^{L_j}s_i(\phi_p,\psi_j^{(i)})\psi_j^{(i)}
 =\Pi_{aux}^{(i)}\tilde\phi_p \qquad\text{for}\quad 1\leq i\leq m
\end{equation*}
 and hence $\phi_p-\tilde\phi_p$ belongs to $\KER$ defined in \eqref{eq:Kernel}.
\par
 Since the functions
 $\phi_p-\tilde\phi_p$ satisfy the Kronecker  property
 $$(\phi_p-\tilde\phi_p)(p')=\begin{cases}
  1&\qquad\text{if $p=p'$}\\[4pt]
  0&\qquad\text{if $p\neq p'$}
 \end{cases}\qquad\text{for $p,p'\in \cV_h\setminus\cN_{aux}$},
 $$
they are linearly independent and the total number of these functions is precisely
$\text{dim}\,V_h-L=\text{dim}\,\KER$ by \eqref{eq:KernelDimension}.
\end{proof}
\begin{remark}\label{rem:KernelBasis}\rm
 The function $\tilde\phi_p$ is supported (i) on one element of $\cT_H$ if $p$ is
 interior to the element,
 (ii) on two elements of $\cT_H$ if $p$ is in the interior of the common edge of
  the two elements, or (iii)
 on four elements of $\cT_H$ if
 $p$ is the common vertex of these elements.
 These three possibilities are illustrated by the blue dots in Figure~\ref{fig:Locations}.
\end{remark}

%
\subsection{A Basis for $\MS$}\label{subsec:MS}
 We are now ready to construct a basis of $\MS$.
 Let the projection $\fC_h:V_h\longrightarrow\KER$ be defined by
\begin{equation}\label{eq:Corrector}
  a(\fC_h v,w)=a(v,w) \qquad \forall\,w\in \KER.
\end{equation}
\begin{remark}\label{rem:Projectionofu}\rm
  Comparing \eqref{eq:uh}, \eqref{eq:IdealMS}, \eqref{eq:IdealMsFEM} and
  \eqref{eq:Corrector}, we see that
\begin{equation}\label{eq:MSError}
  u_h-\MSu=\fC_h u_h.
\end{equation}
\end{remark}
\par
 We can use the normalized hat functions associated with the dual nodes
 (cf. Section~\ref{subsec:Ideal}) and the correction operator $\fC_h$ to construct a basis for $\MS$.
\begin{lemma}\label{lem:IdealBasis} The set
  $\{\hphi_j^{(i)}-\fC_h\hphi_j^{(i)}:\,1\leq i\leq m,1\leq j\leq L_i\}$ is
  a basis of
 $\MS$, where $\hphi_j^{(i)}$ is the normalized hat function corresponding to
  a dual node $\hp_j^{(i)}$ in
 the interior of $K_i$.
\end{lemma}
\begin{proof} It follows from \eqref{eq:Corrector} that $\hphi_j^{(i)}-\fC_h\hphi_j^{(i)}$
 belongs to
$\MS$ defined by \eqref{eq:IdealMS}.
\par
 Note that the functions in
 $\{\hphi_j^{(i)}-\fC_h\hphi_j^{(i)}:\,1\leq i\leq m,1\leq j\leq L_i\}$
are linearly independent because we have the relations
 (cf. \eqref{eq:SameSpan} and \eqref{eq:DualFunctionProjecion})
\begin{align*}
  &\,\text{span}\{\Pi_{aux}^{(i)}(\hphi_j^{(i)}-\fC_h\hphi_j^{(i)}):
  \,1\leq i\leq m,1\leq j\leq L_i\}\\
  =&\,\text{span}\{\Pi_{aux}^{(i)}\hphi_j^{(i)}:\,1\leq i\leq m,1\leq j\leq L_i\}\\
  =&\,\text{span}\{\Pi_{aux}^{(i)}\tilde\phi_j^{(i)}:\,1\leq i\leq m,1\leq j\leq L_i\}\\
  =&\,\text{span}\{\psi_j^{(i)}:\,1\leq i\leq m, 1\leq j\leq L_i\}
\end{align*}
that imply (cf. \eqref{eq:ExplicitDimensionMS})
$$\text{dim}\big(\text{span}\{\hphi_j^{(i)}-
\fC_h\hphi_j^{(i)}:\,1\leq i\leq m,1\leq j\leq L_i\}\big)
=\sum_{i=1}^m L_i=\text{dim}\,\MS.$$

\end{proof}
%
\section{A localized  multiscale finite element method}\label{sec:LOD}
 We will use a special basis of $\KER$ in the construction of the
 localized multiscale finite element space.
 Below $\bA\in\R^{n\times n}$ is the symmetric positive definite
  stiffness matrix representing the
 bilinear form
 $a(\cdot,\cdot)$ with respect to the standard nodal basis functions of $V_h$,
  where $n=\mathrm{dim} V_h$.
%
\subsection{The Modified Gram-Schmidt Algorithm}\label{subsec:MGS}
 A set of vectors $\bw_1,\ldots,\bw_q$ are $\bA$-orthogonal if
 $\bw_j^T\bA\bw_k=0$ for $j\neq k$.  They are $\bA$-orthonormal if in addition they
 are $\bA$-unit vectors, i.e., $\bw_j^T\bA\bw_j=1$ for $j=1,\ldots,q$.
\par
 Given linearly independent vectors $\bv_1,\ldots,\bv_q\in\R^n$, the
 modified Gram-Schmidt algorithm (cf. Algorithm~\ref{alg:GS})
  overwrites $\bv_1,\ldots,\bv_q$ so that the resulting vectors are
  $\bA$-orthonormal and have the same span as the original $\bv_1,\ldots,\bv_q$.
\begin{algorithm}[htbp]
\caption{The modified Gram-Schmidt algorithm.}
\label{alg:GS}
\begin{align*}
&\hspace{-40pt}\text{for $i=1$ to $q$}\\
 r_{ii}&=\sqrt{\bv_i^T\bA\bv_i}\hspace{266pt}\\
    \bv_i&\leftarrow \bv_i/r_{ii}\\
&\hspace{-12pt}\text{for $j=i+1$ to $q$}\\
   &\hspace{10pt}r_{ij}=\bv_i^T\bA\bv_j\\
  &\hspace{10pt} \bv_j\leftarrow \bv_j-\bv_ir_{ij}\\
 &\hspace{-10pt}\text{end}\\
&\hspace{-40pt}\text{end}
\end{align*}
\end{algorithm}
\par
 Properties of the modified Gram-Schmidt algorithm can be found in
 \cite{RTSK:2012:Nonstandard}.
\subsection{Construction of the Special Basis Functions of $\KER$}\label{subsec:Basis}
 We apply a domain decomposition substructuring  strategy to construct three
 types of basis vectors/functions
  according to the three types of functions $\tilde\phi_p$
 in Remark~\ref{rem:KernelBasis}.
%
\subsubsection{Basis Vectors/Functions Associated with the Interior of an Element of $\cT_H$}
\label{subsubsec:Element}
  Let  the
  functions $\phi_{p_1}-\tilde\phi_{p_1},\ldots,\phi_{p_{q_i}}-\tilde\phi_{p_{q_i}}$ be represented by
  the vectors $\bv_1,\ldots,\bv_{q_i}\in\R^n$ with respect to the standard nodal basis of $V_h$, where
  $p_1,\ldots,p_{q_i}$ are the vertices
   of $\cT_h$ interior to
  $K_i\in\cT_H$ and not a dual node (cf. Section~\ref{subsec:Ideal}).  We apply
  the modified Gram-Schmidt  algorithm in Section~\ref{subsec:MGS} to
  replace $\bv_1,\ldots,\bv_{q_i}$ by
  a set of $\bA$-orthonormal vectors.
  The corresponding finite element functions are supported on $K_i$,
  just like the original basis functions
  $\phi_{p_1}-\tilde\phi_{p_1},\ldots,\phi_{p_{q_i}}-\tilde\phi_{p_{q_i}}$.
\begin{remark}\label{rem:CostElement}\rm
  The computational cost of this process is $O((n/m)^3)$ for each of the
   $m$ elements of $\cT_H$ and these
  computations can be carried out in parallel.
\end{remark}
%
\subsubsection{Basis Vectors/Functions Associated with the Interior of an Edge of $\cT_H$ }
\label{subsubsec:Edge}
 Let $E$ be an interior edge of  $\cT_H$ shared by
 $K_1^E,K_2^E\in\cT_H$.  Given a vertex $p$ of $\cT_h$ in the interior of
 $E$, we first replace the vector representing $\phi_p-\tilde\phi_p$ by a vector
 that is $\bA$-orthogonal to  the vectors constructed in
 Section~\ref{subsubsec:Element} for the elements $K_1^E$ and $K_2^E$.
  Then we apply the modified Gram-Schmidt algorithm
  in Section~\ref{subsec:MGS} to replace these vectors by $n_E$ $\bA$-orthonormal
  vectors, where $n_E$ is the number of vertices of $\cT_h$ interior to $E$.
  The corresponding finite element functions are supported on $K_1^E\cup K_2^E$,
   just like the original basis functions.
\begin{remark}\label{rem:CostEdge}\rm
  The computational cost of this process is $O((n/m)^{5/2})$ for each of the $O(m)$
  many interior edges of $\cT_H$ and these
  computations can be carried out in parallel.
\end{remark}
\subsubsection{Basis Vectors/Functions Associated with an Interior Vertex of $\cT_H$}
\label{subsubsec:Vertex}
 Let $p$ be an interior vertex of $\cT_H$ that is the common vertex of the four elements
 $K_1^p,\ldots,K_4^p$ of $\cT_H$
 and the four edges
 $E_1^p,\dots,E_4^p$ of $\cT_H$.  We replace the vector representing $\phi_p-\tilde\phi_p$
 by a $\bA$-unit vector that is $\bA$-orthogonal to the vectors constructed in
 Section~\ref{subsubsec:Element}
 for $K_1^p,\ldots,K_4^p$ and the vectors constructed in
 Section~\ref{subsubsec:Edge} for $E_1^p,\dots,E_4^p$.  The corresponding finite element function is
 supported on $K_1^p\cup K_2^p\cup K_3^p\cup K_4^p$, just like the original basis
 function $\phi_p-\tilde\phi_p$.
\begin{remark}\label{rem:CostVertex}\rm
  The computational cost of this process is $O((n/m)^2)$ for each of the $O(m)$ many
  interior vertices of
  $\cT_H$ and these
  computations can be carried out in parallel.
\end{remark}
\begin{remark}\label{rem:Support}\rm
  By construction each function $\phi_p-\tilde\phi_p$ (where $p$ is not a dual node)
  is replaced by a normalized finite element function supported in the same elements
  from $\cT_H$ as
  $\phi_p-\tilde\phi_p$.
\end{remark}
\subsection{A Matrix Form of the Correction Equation \eqref{eq:Corrector}}
\label{subsec:Correction}
%
 Let $\ell=n-L=\mathrm{dim}\,\KER$ and $\bK\in\R^{n\times \ell}$ be the matrix with full rank
  whose columns contain the
 basis vectors constructed in Section~\ref{subsec:Basis}.  We can then write the correction equation
 \eqref{eq:Corrector} as
\begin{equation*}
   \mathbf{K}^T\bA\bK \bx=\bK^T\bA\bv,
\end{equation*}
 where $\bK\bx$ is the vector representing $\fC_h v$ with respect to the standard nodal basis of $V_h$
 and $\bv$ is the vector representing $v\in V_h$.
 By construction the matrix $\bK^T\bA\bK\in\R^{\ell\times\ell}$ is symmetric positive definite.
 \par
 Since all the columns of $\bK$ are $\bA$-unit vectors by the construction in Section~\ref{subsec:Basis},
 the diagonal part of $\bK^T\bA\bK$ is the identity matrix $\bI\in\R^{\ell\times\ell}$, i.e.,
 \begin{equation}\label{eq:Sructure}
   \bK^T\bA\bK=\bI+\bB,
 \end{equation}
  where $\bB\in\R^{\ell\times\ell}$ is a symmetric off-diagonal matrix.
  Moreover $\bB(i,j)=0$  except
  (i) $\bK(:,i)$ is a basis vector associated with a vertex of $\cT_h$ interior to an edge
  $E_i$  of $\cT_H$
  constructed in Section~\ref{subsubsec:Edge}, $\bK(:,j)$ is a basis vector associated with
  a vertex of $\cT_h$ interior to another edge  $E_j$ of $\cT_H$, and $E_i$ and $E_j$ are
  both edges of one of the elements in $\cT_H$; (ii) $\bK(:,i)$ is a basis vector associated with
  an interior vertex $p_i$ of $\cT_H$ constructed in Section~\ref{subsubsec:Vertex},
  $\bK(:,j)$ is a basis vector associated with another  interior vertex $p_j$ of $\cT_H$, and $p_i$
  and $p_j$ are endpoints of  an edge $E$ of  $\cT_H$.
 Therefore $\bK^T\bA\bK$ is a sparse matrix because of \eqref{eq:Sructure}.
\par
 We will denote the condition number of $\bK^T\bA\bK$ by $\mathfrak{K}$, i.e,
\begin{equation}\label{eq:ConditonNUmberBdd}
  \mathfrak{K}=\frac{\lambda_{\max}(\bK^T\bA\bK)}{\lambda_{\min}(\bK^T\bA\bK)}.
\end{equation}
\begin{remark}\label{rem:ConditionNumber}\rm
     $\mathfrak{K}$ can be estimated through conjugate gradient iterations.
\end{remark}
\begin{remark}\label{rem:Preconditioning}\rm
  Numerical results indicate that the condition number $\mathfrak{K}$ is moderate even for
  large $\beta$ and it
  decreases as $H$ decreases.  Heuristically we can treat $\bK^T\bA\bK$ as a preconditioned
   matrix where
  the preconditioner is similar to a substructuring domain decomposition preconditioner
  (cf. \cite{BPS:1986:I})
  where global communication is implied by the definition of $\KER$.
\end{remark}
\subsection{A Localized Multiscale Finite Element Method}\label{subsec:LOD}
 Let $\mathbf{f}_j^{(i)}\in\R^n$ ($1\leq i\leq m, 1\leq j\leq L_i$) be the vector that represents the
 normalized hat
 function $\hphi_j^{(i)}$ associated with the dual node $\hp_j^{(i)}$ (cf. Section~\ref{subsec:Ideal})
 with respect to the standard nodal basis of $V_h$.
 We apply $k$ iterations of the conjugate gradient algorithm
 (cf. \cite{Greenbaum:1997:Iterative,Saad:2003:IM}) with initial guess $\mathbf{0}$ to the corrector equation
\begin{equation}\label{eq:CorectorEq}
  \bK^T\bA\bK\, \bx =\bK^T\bA\, \mathbf{f}_j^{(i)}
\end{equation}
 to obtain an approximate solution $\hat\bx_j^{(i)}$ and define $\fC_{h,k}\hphi_j^{(i)}$
 to be the finite element function represented by $\bK\hat\bx_j^{(i)}$.  We can then
  extend $\fC_{h,k}$ to be a linear operator from $\Vd$ to $\KER$ by linearity.
\begin{remark}\label{rem:PCG}\rm
  We can regard $\bK\hat\bx_j$ as an approximate solution for $\fC_h v$
  in \eqref{eq:Corrector} (with $v=\hphi_j^{(i)}$)
   obtained by a preconditioned conjugate gradient
  algorithm.
\end{remark}
\begin{remark}\label{rem:LOD}\rm
  According to Remark~\ref{rem:Support},
   the function $\fC_{h,k}\hphi_j^{(i)}$ is supported on a domain obtained by adding
   roughly $k$ layers of
  elements in $\cT_H$ to $K_i$, i.e., it is a localized correction of $\hphi_j^{(i)}$.
\end{remark}
\par
 Let $\MSk$ be the space of functions spanned by $\hphi_j^{(i)}-\fC_{h,k}\hphi_j^{(i)}$
 for $1\leq i\leq m$ and $1\leq j\leq L_i$.  The localized  multiscale finite element method for
 \eqref{eq:BVP} is to find $\MSku\in\MSk$ such that
\begin{equation}\label{eq:LOD}
  a(\MSku,v)=\int_\O fv\,dx \qquad\forall\,v\in\MSk.
\end{equation}
\begin{remark}\label{rem:MSk}\rm
 Since $\Pi_{aux}\big(\hphi_j^{(i)}-\fC_{h,k}\hphi_j^{(i)}\big)=\Pi_{aux}\hphi_j^{(i)}$, the functions
 $\hphi_j^{(i)}-\fC_{h,k}\hphi_j^{(i)}$
 for $1\leq i\leq m$ and $1\leq j\leq L_i$ are linearly independent by Remark~\ref{rem:DualFunctions}
 and hence $\d\mathrm{dim}\MSk=\mathrm{dim}\MS$.
\end{remark}
\begin{remark}\label{rem:ROM}\rm
  The localized multiscale finite element method is a reduced order method. The time
   consuming construction
  of the basis of the space $\MSk$ can be carried out in the off-line stage.
  The on-line solution
  of \eqref{eq:LOD} is very fast when the dimension of $\MSk$ is moderate.
\end{remark}
\section{Analysis of the Localized Multiscale Finite Element Method}
\label{sec:Analysis}
 It follows from \eqref{eq:NormalizedHatFunction}, \eqref{eq:Corrector}, \eqref{eq:CorectorEq} and
 the theory of the conjugate gradient algorithm (cf. \cite{Greenbaum:1997:Iterative,Saad:2003:IM}) that
\begin{equation}\label{eq:CGEst2}
  \|\fC_h\hphi_j^{(i)}-\fC_{h,k}\hphi_j^{(i)}\|_a\leq \frac{2q^k}{1+q^{2k}}\|\fC_h\hphi_j^{(i)}\|_a\leq
   \frac{2q^k}{1+q^{2k}}\|\hphi_j^{(i)}\|_a=\frac{2q^k}{1+q^{2k}},
\end{equation}
 where, in view of \eqref{eq:ConditonNUmberBdd}, we can take
\begin{equation}\label{eq:qDef}
  q=\frac{\sqrt{\mathfrak{K}}-1}{\sqrt{\mathfrak{K}}+1}.
\end{equation}
%
 %
 \subsection{Estimate for $\|\fC_h-\fC_{h,k}\|_a$}\label{subsec:fChk}
 Let $\|\fC_h-\fC_{h,k}\|_a$  be the norm induced by the energy norm for the operator
 $\fC_h-\fC_{h,k}$
 that maps $\Vd$ to $\KER$.
\begin{lemma}\label{lem:fChk}
   We have the estimate
\begin{equation}\label{eq:DifferenceBdd}
  \|\fC_h-\fC_{h,k}\|_a\leq \frac{2q^k}{1+q^{2k}} \sqrt{L}.
\end{equation}
\end{lemma}
\begin{proof}  It follows from \eqref{eq:DimensionVaux}, \eqref{eq:DualKronecker},
 \eqref{eq:CGEst2} and
the Cauchy-Schwarz inequality that
\begin{align*}
  \big\|(\fC_h-\fC_{h,k})\sum_{i=1}^m\sum_{j=1}^{L_i}c_j^{(i)}\hphi_j^{(i)}\big\|_a
  &\leq \sum_{i=1}^m\sum_{j=1}^{L_i}|c_j^{(i)}|\,\|(\fC_h-\fC_{h,k})\hphi_j^{(i)}\|_a\\
  &\leq \frac{2q^k}{1+q^{2k}}\sum_{i=1}^m\sum_{j=1}^{L_i}|c_j^{(i)}|\\
  &\leq \frac{2q^k}{1+q^{2k}} \Big(\sum_{i=1}^m L_i\Big)^\frac12
  \Big(\sum_{i=1}^m\sum_{j=1}^{L_i}|c_j^{(i)}|^2\Big)^\frac12\\
  &= \frac{2q^k}{1+q^{2k}}\sqrt{L}\big\|\sum_{i=1}^m\sum_{j=1}^{L_i}c_j^{(i)}\hphi_j^{(i)}\big\|_a,
\end{align*}
 which implies \eqref{eq:DifferenceBdd}.
\end{proof}
%
\subsection{Error Estimates for $\MSku$}\label{subsec:MSku}
 We will denote by $M$ the number $\d\max_{1\leq i\leq m} M_i$, where $M_i$ is
 defined in \eqref{eq:Mi}.
\begin{theorem}\label{thm:AbstractEnergyError}
  The solution $\MSku$ of \eqref{eq:LOD} satisfies the estimate
\begin{equation}\label{eq:AbstractEnergyError}
  \|u_h-\MSku\|_a\leq \|u_h-\MSu\|_a+\frac{2q^k}{1+q^{2k}}
  \sqrt{L}\sqrt{M} H^{-1}\sqrt{\beta} \|f\|_\LT,
\end{equation}
 where $u_h\in V_h$ $($resp., $\MSu)$ satisfies \eqref{eq:uh}
 $($resp., \eqref{eq:IdealMsFEM}$)$, and
 $q$ is defined in \eqref{eq:qDef}.
\begin{proof}
  Let $u_h^{(d)}\in V_h^d$ be defined by
\begin{equation}\label{eq:uhd}
  u_h^{(d)}=\sum_{i=1}^m\sum_{j=1}^{L_i}s_i(u_h,\psi_j^{(i)})\tilde\phi_j^{(i)}.
\end{equation}
Then we have, in view of \eqref{eq:LocalPiAux}, \eqref{eq:PiAuxVh} and
\eqref{eq:DualFunctionProjecion},
\begin{equation*}
  \Pi_{aux} u_h^{(d)}=\Pi_{aux} u_h,
\end{equation*}
 and hence $u_h-u_h^{(d)}$ belongs to $\KER$ so that
\begin{equation}\label{eq:Trivial}
  \fC_h (u_h-u_h^{(d)})=u_h-u_h^{(d)}.
\end{equation}
\par
 Note that by construction (cf. Section~\ref{subsec:LOD})
\begin{equation}\label{eq:Membership}
   u_h^{(d)}-\fC_{h,k} u_h^{(d)}\quad\text{belongs to}\quad \MSL.
\end{equation}
\par
 It follows from \eqref{eq:MSError},
 \eqref{eq:Trivial},  \eqref{eq:Membership}, the Galerkin orthogonality of the LOD multiscale
  finite element method and the triangle inequality that
\begin{align}\label{eq:Abstract1}
  \|u_h-\MSLu\|_a&\leq \|u_h-(u_h^{(d)}-\fC_{h,k} u_h^{(d)})\|_a\\
                 &\leq \|\fC_h u_h\|_a+
                     \|(\fC_h-\fC_{h,k})u_h^{(d)}\|_a
                     =\|u_h-\MSu\|_a+\|(\fC_h-\fC_{h,k})u_h^{(d)}\|_a,\notag
\end{align}
and we have
\begin{equation}\label{eq:Abstract2}
  \|(\fC_h-\fC_{h,k})u_h^{(d)}\|_a \leq \frac{2q^k}{1+q^{2k}}
  \sqrt{L}\|u_h^{(d)}\|_a
\end{equation}
 by \eqref{eq:DifferenceBdd}.
\par
{
 From \eqref{eq:Comparison},
  \eqref{eq:PF}, \eqref{eq:uhBdd}, \eqref{eq:LocalLTwo}, \eqref{eq:Mi}
   and \eqref{eq:uhd} we find
\begin{align*}
 \|u_h^{(d)}\|_a^2&=\sum_{i=1}^m a\Big(\sum_{j=1}^{L_i}s_i(u_h,\psi_j^{(i)})\tilde\phi_i^{(j)},
  \sum_{\ell=1}^{L_i}s_i(u_h,\psi_\ell^{(i)})\tilde\phi_\ell^{(j)}\Big)\\
  &\leq \sum_{i=1}^m M_i\sum_{j=1}^{L_i}[s_i(u_h,\psi_j^{(i)})]^2\\
  &\leq M\sum_{i=1}^m\sum_{j=1}^{L_i}[s_i(u_h,\psi_j^{(i)})]^2\\
  &\leq M H^{-2}\|u_h\|_{L^2(\O;\kappa)}^2\\
  &\leq M H^{-2}\beta\|u_h\|_\LT^2\\
  &\leq  M H^{-2}\beta \|u_h\|_a^2\\
   &\leq  M H^{-2}\beta \|f\|_\LT^2,
\end{align*}
 and hence
\begin{equation}\label{eq:Abstract3}
  \|u_h^{(d)}\|_a\leq \sqrt{M} H^{-1}\sqrt\beta\|f\|_\LT.
\end{equation}
}
\par
 The estimate \eqref{eq:AbstractEnergyError} follows from
  \eqref{eq:Abstract1}--\eqref{eq:Abstract3}.
\end{proof}
\end{theorem}
\par
 Putting Theorem~\ref{thm:IdealErrors} and Theorem~\ref{thm:AbstractEnergyError}
 together, we arrive at the
 estimate
\begin{equation}\label{eq:ExplicitEneryError}
  \|u_h-\MSku\|_a\leq \Big[C_*H+\frac{2q^k}{1+q^{2k}}
  \sqrt{L}\sqrt{M} H^{-1}\sqrt{\beta}\Big]\|f\|_\LT,
\end{equation}
 where the constants $L$, $C_*$, $M$ and $q$ are available through  the construction of
 $V_{aux}^{(i)}$, \eqref{eq:C*}, \eqref{eq:Mi},
  Remark~\ref{rem:ConditionNumber}, and \eqref{eq:qDef}.
\par
 It follows from \eqref{eq:ExplicitEneryError} that
\begin{equation}\label{eq:ConcreteEnergyError}
  \|u_h-\MSLu\|_a\leq (C_*+1) H \|f\|_\LT
\end{equation}
 if $k$ is sufficiently large so that
\begin{equation}\label{eq:ChoiceOfk}
  2q^k\sqrt{L}\sqrt{M}\sqrt{\beta} \leq H^2.
\end{equation}
\par
 A standard duality argument then yields
\begin{equation}\label{eq:ConcreteLTError}
   \|u_h-\MSLu\|_{\LT}\leq [(C_*+1) H]^2\|f\|_\LT.
\end{equation}
\begin{remark}\label{rem:Standard}\rm
  The estimates \eqref{eq:ConcreteEnergyError} and \eqref{eq:ConcreteLTError}
  indicate that the performance
  of the localized multiscale finite element method is similar to the performance
  of standard finite element methods for the
  {homogeneous Dirichlet boundary value problem for the
  Poisson equation on smooth or convex domains.}
\end{remark}
\section{Numerical Results}\label{sec:Numerics}
 The domain $\O$ in the numerical examples is the unit square $(0,1)\times(0,1)$
  and
 the function $f(x)$ in all the examples is defined by
\begin{equation}\label{eq:fDef}
f(x)=\begin{cases}
  0&\qquad\text{for}\;x\in[0,\frac12)\times[0,1]\\[4pt]
  1&\qquad\text{for}\;x\in[\frac12,1]\times[0,1]
\end{cases}
\end{equation}
 and $\|f\|_\LT=1/2$.
\par
 For each example we obtain the following data from the computation:
\begin{itemize}\itemsep=5pt
\item The number $\sqrt{L}$, where $\d L=\sum_{i=1}^m L_i$ and $L_i$ is the number of
eigenfunctions in the local auxiliary space $V_{aux}^{(i)}$ (cf. Section~\ref{sec:Auxiliary}).
\item The number $\sqrt{M}$, where $\d M=\max_{1\leq i\leq m} M_i$ and $M_i$ is defined
 in \eqref{eq:Mi}.
\item The contraction number $q$ of the conjugate gradient iteration defined in \eqref{eq:qDef}.
\item The smallest integer $k$  that satisfies \eqref{eq:ChoiceOfk}.
\end{itemize}
 Then we construct the basis functions of the localized  multiscale finite element space
 $\MSk$ {(cf. Section~5.1--Section~5.3)}.
  {Finally for each $K_i$ we run the modified Gram-Schmidt algorithm (Algorithm~5.1) for the basis functions
 of $\MSk$ that are associated with the dual nodes in each $K_i$, which generates a new set of basis
 functions of  $\MSk$.}
All these can be
 carried out in the off-line stage.
\par
In the on-line stage, the equation \eqref{eq:LOD} is solved  by using the
 basis of $\MSk$ generated at the end of the off-line stage.  {The solution by
  the backslash command in MATLAB
 is very fast (cf. Table~\ref{tab:PSTimes}, Table~\ref{tab:IrregularTimes} and
 Table~\ref{tab:GFTimes}).}
\par
 The computations were carried out on a Dell Inspiron 13 Laptop (11th Gen Intel(R) Core(TM)
  i5-11320H @ 3.20GHz (2.50 GHz), RAM: 8.00GB).  The conjugate gradient iteration used a tolerance of
  $10^{-14}$ for the relative error as the stopping criterion.
\begin{example}\label{example:PS}
\rm
 This is the example from \cite[Section~5.2]{PS:2016:Contrast} with four high contrast channels
 in the unit square (cf. Figure~\ref{fig:PS}).
\begin{figure}[htbp]
\begin{center}
  \includegraphics[height=2in]{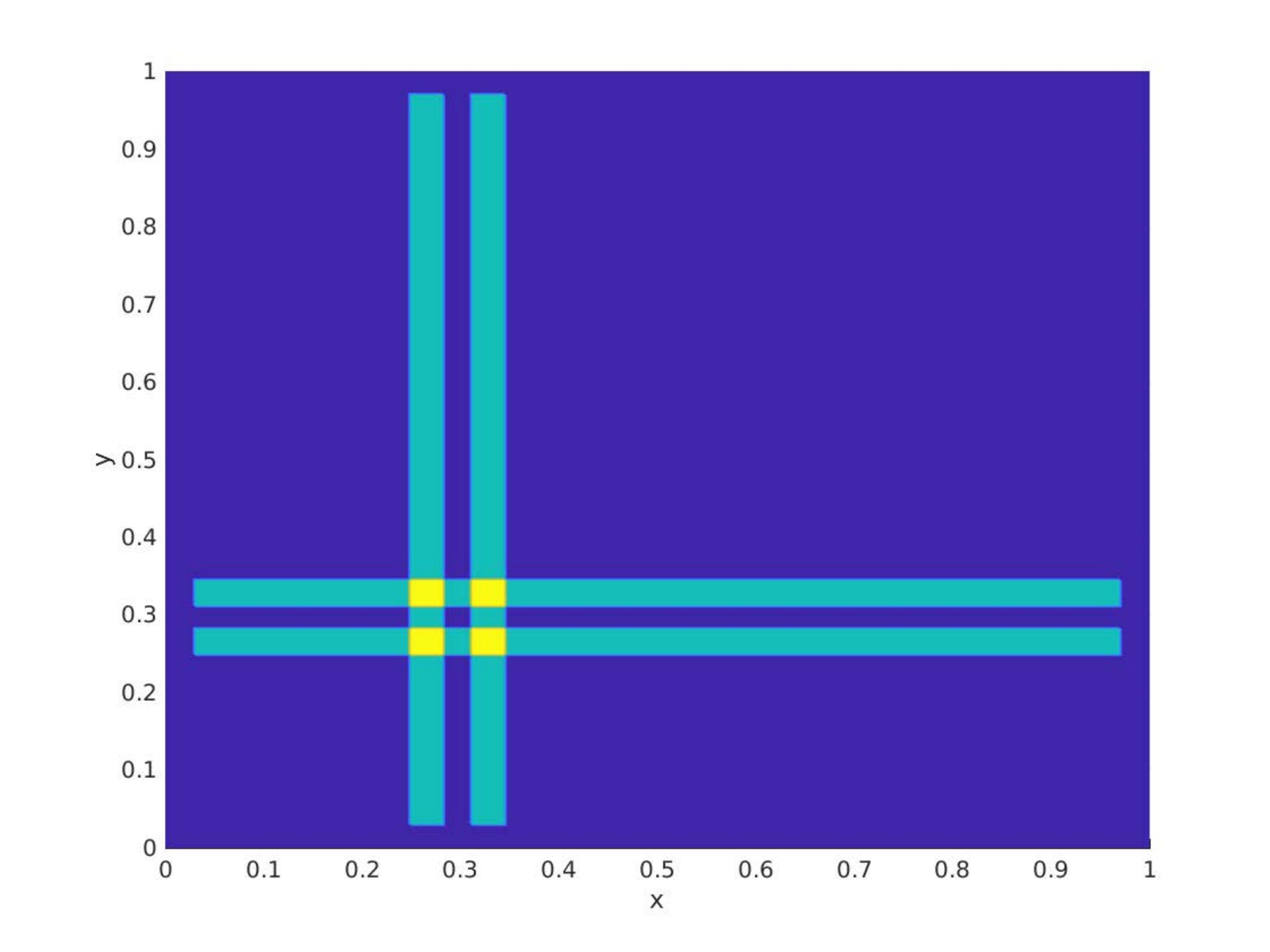}
\end{center}
\par\vspace{-.2in}
\caption{Four high contrast channels.}
\label{fig:PS}
\end{figure}
\par
 The diffusion coefficient $\kappa(x)$ is defined by
 $$\kappa(x)=A(x_1,x_2)+A(x_2,x_1),$$
 where
 $$A(x)=\begin{cases}
   \beta/2 &\qquad\text{for}\;x\in[\frac{8}{32},\frac{9}{32}]\times [\frac{1}{32},\frac{31}{32}]
   \cup[\frac{10}{32},\frac{11}{32}]\times[\frac{1}{32},\frac{31}{32}]\\[4pt]
   1 &\qquad\text{elsewhere}
 \end{cases}\,.
 $$
\par
 We solve the boundary value problem \eqref{eq:BVP} for $\beta=10^2,10^4,10^6$ and $10^8$.  We take
 $h=1/256$ and $H=1/8, 1/16, 1/32, 1/64$.
 The data generated by the computation  are displayed in Table~\ref{tab:PSData}.
{
\begin{table}[htbp]
\resizebox{1 \textwidth}{!}{
  \begin{tabular}{|l||c|c|c|c||c|c|c|c||c|c|c|c||c|c|c|c|}
  \hline
  \multirow{2}{*}{\backslashbox{$H$}{$\beta$}}&\multicolumn{4}{c||}{\mystrut(12,4)$10^2$}&\multicolumn{4}{c||}{$10^4$}&
  \multicolumn{4}{c||}{$10^6$}& \multicolumn{4}{c|}{$10^8$}\\
  \cline{2-17}&\mystrut(14,4)$\sqrt{L}$&$\sqrt{M}$&$q$&$k$&$\sqrt{L}$&$\sqrt{M}$&$q$&$k$&$\sqrt{L}$&$\sqrt{M}$&$q$&$k$
  &$\sqrt{L}$&$\sqrt{M}$&$q$&$k$\\
   \hline
   &&&&&&&&&&&&&&&&\\[-9pt]
  $1/8$ &$8.6$&$2.1\times10^6$&$0.54$&$39$&$8.6$&$1.6\times10^5$&$0.56$&$41$&$8.6$&$9.9\times10^{5}$&$0.56$&$48$&
  $8.6$&$1.0\times10^7$&$0.56$&$56$\\
  \hline &&&&&&&&&&&&&&&&\\[-10pt]
  $1/16$ &$16$&$4.2\times10^3$&$0.44$&$24$&$24$&$4.2\times10^4$&$0.44$&$30$&$16$&$4.1\times10^5$&$0.44$&$36$
  &$16$&$4.1\times10^6$&$0.44$&$42$\\
  \hline&&&&&&&&&&&&&&&&\\[-10pt]
  $1/32$&$32$&$2.5\times10^2$&$0.40$&$21$&$32$&$2.5\times10^2$&$0.40$&$24$&$32$&$2.5\times10^2$&$0.40$&$26$
  &$32$&$2.5\times10^2$&$0.40$&$29$\\
  \hline&&&&&&&&&&&&&&&&\\[-10pt]
  $1/64$&$64$&$5.0\times10^1$&$0.38$&$21$&$64$&$5.0\times10^1$&$0.38$&$23$&$64$&$5.0\times10^1$&$0.38$&$26$
  &$64$&$5.0\times10^1$&$0.38$&$28$\\
  \hline
  \end{tabular}
  }
\par\vspace{.1in}
\caption{Data generated by the computation for Example~\ref{example:PS}.}
\label{tab:PSData}
\end{table}
}
\par
 The estimates for the energy error and $L^2$ error according to
  \eqref{eq:ConcreteEnergyError} and
 \eqref{eq:ConcreteLTError} are shown
 in Table~\ref{tab:PSEstimates},
 and the
  actual energy errors and  $L^2$ errors can be found in Table~\ref{tab:PSErrors}.
  It is observed that for a given $H$, the actual errors are independent of the
   contrast $\beta$.  By comparing
  these two tables we see that the estimates \eqref{eq:ConcreteEnergyError} and
 \eqref{eq:ConcreteLTError} are satisfied.
\begin{table}[H]
\begin{center}
{\scriptsize
  \begin{tabular}{|l|c|c|c|c|}
  \hline&&&&\\[-10pt]
    \hspace{52pt}$H$ & $\frac18$ & $\frac{1}{16}$ & $\frac{1}{32}$ & $\frac{1}{64}$\\
    &&&&\\[-10pt]
    \hline&&&&\\[-10pt]
  Energy Error Estimate&$1.2\times10^{-1}$&$5.9\times10^{-2}$&$3.0\times10^{-2}$&$1.5\times10^{-2}$\\
    \hline&&&&\\[-10pt]
   $L^2$ Error Estimate&$1.4\times10^{-2}$&$3.5\times10^{-3}$&$8.8\times10^{-4}$&$2.2\times10^{-4}$\\
    \hline
  \end{tabular}
}
\end{center}\par\vspace{.1in}
\caption{Estimates \eqref{eq:ConcreteEnergyError} and
 \eqref{eq:ConcreteLTError} for Example~\ref{example:PS}.}
 \label{tab:PSEstimates}
\end{table}
\begin{table}[H]
\resizebox{1 \textwidth}{!}{
  \begin{tabular}{|l||c|c||c|c||c|c||c|c|}
     \hline
  \multirow{2}{*}{\backslashbox{$H$}{$\beta$}}&\multicolumn{2}{c||}{\mystrut(12,4)$10^2$}&\multicolumn{2}{c||}{$10^4$}&
  \multicolumn{2}{c||}{$10^6$}& \multicolumn{2}{c|}{$10^8$}\\
  \cline{2-9}&\mystrut(12,4)Energy Error&$L^2$ Error&Energy Error&$L^2$ Error
  &Energy Error&$L^2$ Error&Energy Error&$L^2$ Error\\
    \hline
   &&&&&&&&\\[-9pt]
  $1/8$&$3.1\times10^{-3}$&$4.8\times10^{-5}$  &$3.1\times10^{-3}$&$4.8\times10^{-5}$   &$3.1\times10^{-3}$&$4.8\times10^{-5}$
  &$3.1\times10^{-3}$&$4.8\times10^{-5}$ \\
  \hline&&&&&&&&\\[-10pt]
  $1/16$&$1.7\times10^{-3}$&$1.6\times10^{-5}$ &$1.7\times10^{-3}$&$1.6\times10^{-5}$&$1.7\times10^{-3}$&$1.6\times10^{-5}$
  &$1.7\times10^{-3}$& $1.6\times 10^{-5}$\\
  \hline&&&&&&&&\\[-10pt]
  $1/32$&$3.5\times10^{-4}$&$1.5\times10^{-6}$  &$3.5\times10^{-4}$&$1.5\times10^{-6}$  &$3.5\times10^{-4}$&$1.5\times10^{-6}$
  &$3.5\times10^{-4}$&$1.5\times10^{-6}$ \\
  \hline&&&&&&&&\\[-10pt]
  $1/64$&$1.1\times10^{-4}$&$2.3\times10^{-7}$ &$1.1\times10^{-4}$&$2.3\times10^{-7}$ &$1.1\times10^{-4}$&$2.3\times10^{-7}$
  &$1.1\times10^{-4}$&$2.3\times10^{-7}$ \\
  \hline
  \end{tabular}
}
\par\vspace{.1in}
\caption{Energy and $L^2$ errors for Example~\ref{example:PS}.}
\label{tab:PSErrors}
\end{table}
\par
The convergence histories in the energy norm and the  $L^2$ norm for the
 standard $Q_1$ finite element method
and the localized  multiscale finite element method are displayed
 in Figure~\ref{fig:FourChannelsErrors}.
The advantage of the multiscale finite element method is clearly observed.
 The order of
convergence of the multiscale finite element method in the energy norm is between
$1$ and $2$, and the
order of convergence in the $L^2$ norm is between $2$ and $3$.
The relative errors of our method also
appear to be smaller than the ones in \cite{PS:2016:Contrast}.
\begin{figure}[htbp]
\begin{center}
  \includegraphics[height=2in]{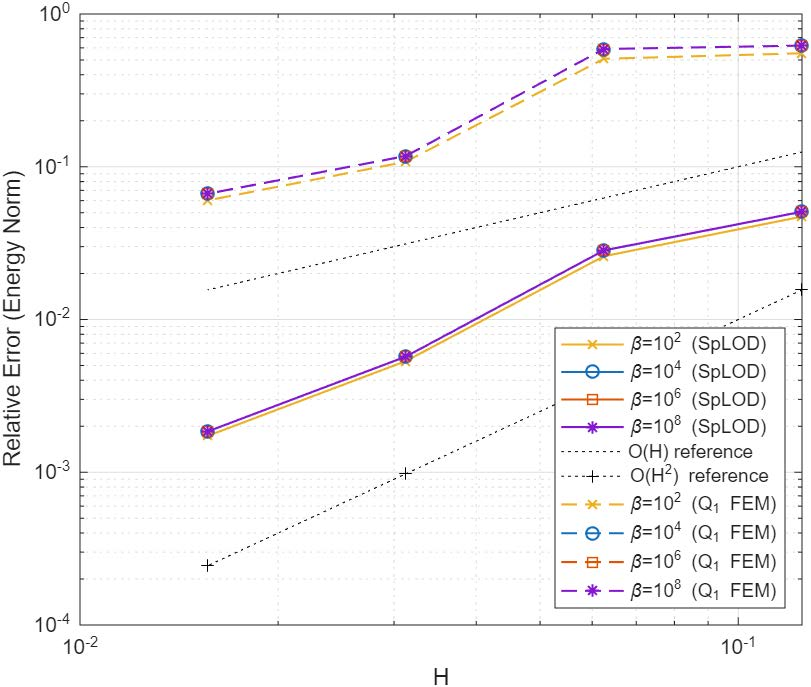} \qquad\qquad
 \includegraphics[height=2in]{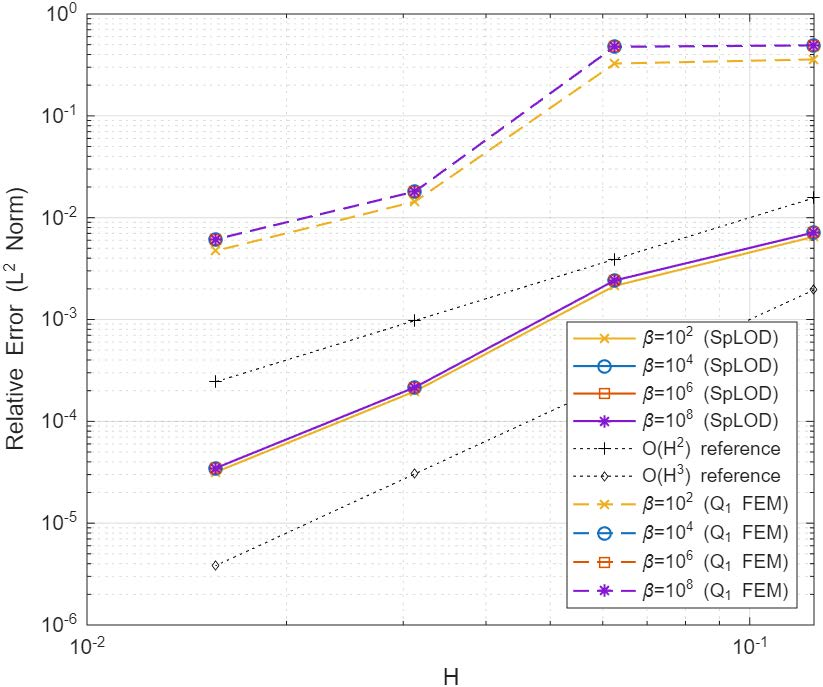}
\end{center}
\caption{Relative energy errors (left) and relative
$L^2$ errors (right) for Example~\ref{example:PS}.}
\label{fig:FourChannelsErrors}
\end{figure}
\end{example}
\par
{
 Finally we report in Table~\ref{tab:PSTimes} the solution times of
  \eqref{eq:LOD} by the backslash command in MATLAB.}
\begin{table}[H]
\footnotesize
\begin{center}
\begin{tabular}{|c|c|c|c|c|}
\hline
&&&&\\[-11pt]
$H \setminus \beta$ &  $10^2$& $10^4$& $10^6$ & $10^8$ \\
\hline
&&&&\\[-11pt]
1/8 &   $7.83\times10^{-5}$ &  $6.95\times10^{-5}$ &  $7.95\times 10^{-5}$&  $1.10\times10^{-4}$\\
\hline
&&&&\\[-11pt]
1/16&   $9.33\times10^{-4}$&  $7.69\times10^{-4}$ &  $1.08\times10^{-3}$&  $9.39\times10^{-4}$\\
\hline
&&&&\\[-11pt]
1/32&   $2.13\times10^{-2}$& $2.14\times10^{-2}$ &  $2.11\times10^{-2}$ &  $1.98\times10^{-2}$\\
\hline
\end{tabular}
\end{center}
\caption{{Solution times (in seconds) of \eqref{eq:LOD} by the backslash command} {in MATLAB
for Example~\ref{example:PS}}.}
\label{tab:PSTimes}
\end{table}

\newpage
\begin{example}\label{example:Irregular}\rm
This is an example of irregular channels in the unit square
 where the diffusion coefficient
$\kappa(x)=\beta$ in the red region and $1$ elsewhere
(cf. Figure~\ref{fig:Irregular}).
\begin{figure}[H]
\begin{center}
  \includegraphics[height=2in]{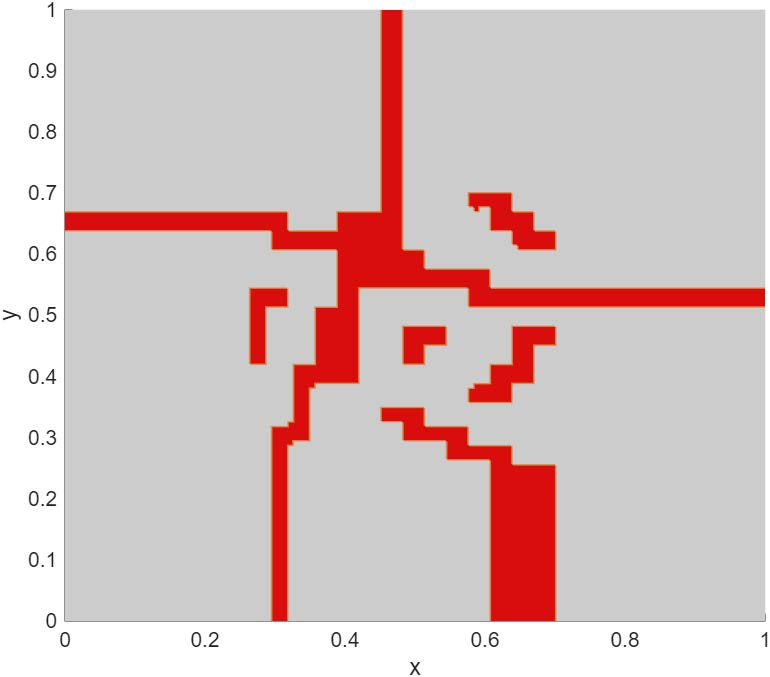}
\end{center}
\par\vspace{-.2in}
\caption{Irregular high contrast channels.}
\label{fig:Irregular}
\end{figure}
\par
We solve the boundary value problem \eqref{eq:BVP} for $\beta=10^2,10^4,10^6$ and $10^8$.  We take
 $h=1/256$ and $H=1/8, 1/16, 1/32, 1/64$.
  The data generated by the computation  are displayed in Table~\ref{tab:IrregularData}.
  \begin{table}[htbp]
  \resizebox{1 \textwidth}{!}{
  \begin{tabular}{|l||c|c|c|c||c|c|c|c||c|c|c|c||c|c|c|c|}
  \hline
  \multirow{2}{*}{\backslashbox{$H$}{$\beta$}}&\multicolumn{4}{c||}{\mystrut(12,4)$10^2$}&\multicolumn{4}{c||}{$10^4$}&
  \multicolumn{4}{c||}{$10^6$}& \multicolumn{4}{c|}{$10^8$}\\
  \cline{2-17}&\mystrut(14,4)$\sqrt{L}$&$\sqrt{M}$&$q$&$k$&$\sqrt{L}$&$\sqrt{M}$&$q$&$k$&
  $\sqrt{L}$&$\sqrt{M}$&$q$&$k$&$\sqrt{L}$&$\sqrt{M}$&$q$&$k$\\
  \hline
   &&&&&&&&&&&&&&&&\\[-9pt]
  $1/8$ &$8.7$&$3.9\times10^5$&$0.46$&$29$&$8.7$&$2.3\times10^6$&$0.46$&$34$&$8.7$&$2.3\times10^7$&$0.46$&$40$&
  $8.7$&$2.3\times10^8$&$0.46$&$46$\\
  \hline
   &&&&&&&&&&&&&&&&\\[-10pt]
  $1/16$ &$17$&$1.5\times10^5$&$0.47$&$31$&$17$&$1.5\times10^5$&$0.47$&$34$&$17$&$6.6\times10^5$&$0.47$&$39$&
  $17$&$6.6\times10^6$&$0.47$&$46$\\
  \hline
   &&&&&&&&&&&&&&&&\\[-10pt]
  $1/32$&$32$&$1.3\times10^4$&$0.41$&$26$&$32$&$2.6\times10^4$&$0.41$&$29$&$32$&$1.4\times10^5$&$0.41$&$34$&
  $32$&$1.4\times10^6$&$0.41$&$39$\\
  \hline
   &&&&&&&&&&&&&&&&\\[-10pt]
  $1/64$&$64$&$1.9\times10^2$&$0.41$&$24$&$64$&$1.8\times10^3$&$0.42$&$29$&$64$&$1.8\times10^4$&$0.42$&$35$&
  $64$&$1.8\times10^5$&$0.42$&$40$\\
  \hline
  \end{tabular}
}
\par\vspace{.1in}
\caption{Data generated by the computation for Example~\ref{example:Irregular}.}
\label{tab:IrregularData}
\end{table}
\par
 The estimates for the energy error and $L^2$ error according to
 \eqref{eq:ConcreteEnergyError} and
 \eqref{eq:ConcreteLTError} are identical to the ones in Table~\ref{tab:PSEstimates},
 and the actual
 energy errors and $L^2$ errors can be found in Table~\ref{tab:IrregularErrors}.
 The behavior of the errors are similar to that in Example~\ref{example:PS},
 namely for a given $H$ the errors are independent of the contrast $\beta$ and
 the estimates \eqref{eq:ConcreteEnergyError} and
 \eqref{eq:ConcreteLTError} are satisfied.
\begin{table}[htbp]
\resizebox{1 \textwidth}{!}{
  \begin{tabular}{|l||c|c||c|c||c|c||c|c|}
     \hline
  \multirow{2}{*}{\backslashbox{$H$}{$\beta$}}&\multicolumn{2}{c||}{\mystrut(12,4)$10^2$}&\multicolumn{2}{c||}{$10^4$}&
  \multicolumn{2}{c||}{$10^6$}& \multicolumn{2}{c|}{$10^8$}\\
  \cline{2-9}&\mystrut(12,4)Energy Error&$L^2$ Error&Energy Error&$L^2$ Error
  &Energy Error&$L^2$ Error&Energy Error&$L^2$ Error\\
    \hline
   &&&&&&&&\\[-9pt]
  $1/8$&$8.6\times10^{-3}$&$2.4\times10^{-4}$ &$8.9\times10^{-3}$&$2.6\times10^{-4}$  &$8.9\times10^{-3}$&$2.6\times10^{-4}$
  &$8.9\times10^{-3}$&$2.6\times10^{-4}$ \\
  \hline&&&&&&&\\[-10pt]
  $1/16$&$2.7\times10^{-3}$&$3.6\times10^{-5}$   &$2.8\times10^{-3}$&$3.8\times10^{-5}$ &$2.8\times10^{-3}$&$3.8\times10^{-5}$
  &$2.8\times10^{-3}$&$3.8\times10^{-5}$  \\
  \hline&&&&&&&&\\[-10pt]
  $1/32$&$1.5\times10^{-3}$&$1.1\times10^{-5}$  &$1.5\times10^{-3}$&$1.2\times10^{-5}$ &$1.5\times10^{-3}$&$1.2\times10^{-5}$
  &$1.5\times10^{-3}$&$1.2\times10^{-5}$ \\
  \hline&&&&&&&&\\[-10pt]
  $1/64$&$1.9\times10^{-4}$&$6.3\times10^{-7}$  &$1.9\times10^{-4}$&$6.4\times10^{-7}$& $1.9\times10^{-4}$&$6.4\times10^{-7}$
  &$1.9\times10^{-4}$&$6.4\times10^{-7}$   \\
  \hline
  \end{tabular}
}
\par\vspace{.1in}
\caption{Energy and  $L^2$ errors for Example~\ref{example:Irregular}.}
\label{tab:IrregularErrors}
\end{table}
\par
 The convergence histories in the energy norm and the $L^2$ norm for the
 standard $Q_1$ finite element method
 and the localized multiscale finite element method are displayed in
  Figure~\ref{fig:IrregularErrors}.
 Again we observe
 the advantage of the multiscale finite element method over the standard
  $Q_1$ finite element method
 and  the order of convergence of the multiscale finite element method
  in the energy (resp., $L^2$ )
 norm is between $1$ and $2$ (resp., $2$ and $3$).
 \begin{figure}[H]
\begin{center}
  \includegraphics[height=2in]{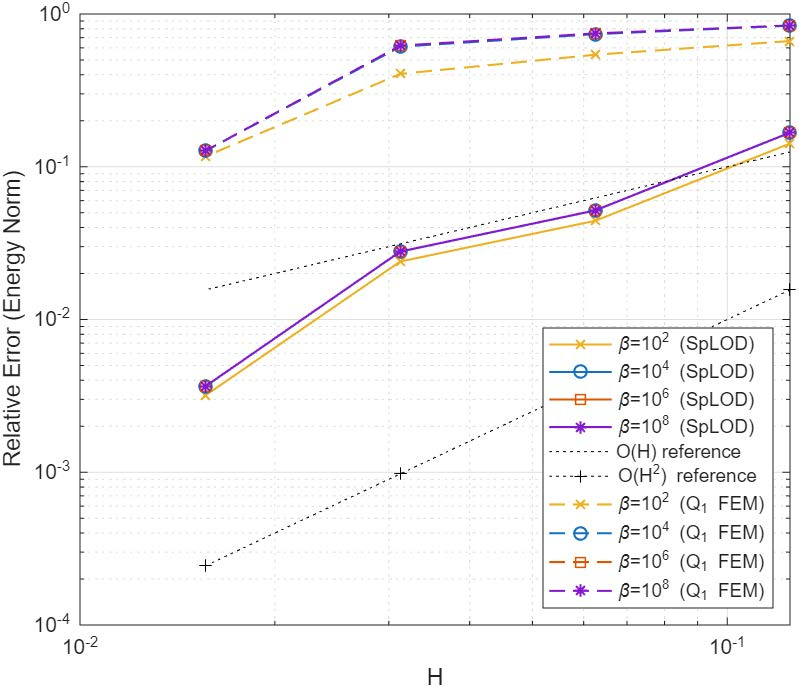} \qquad \qquad \includegraphics[height=2in]{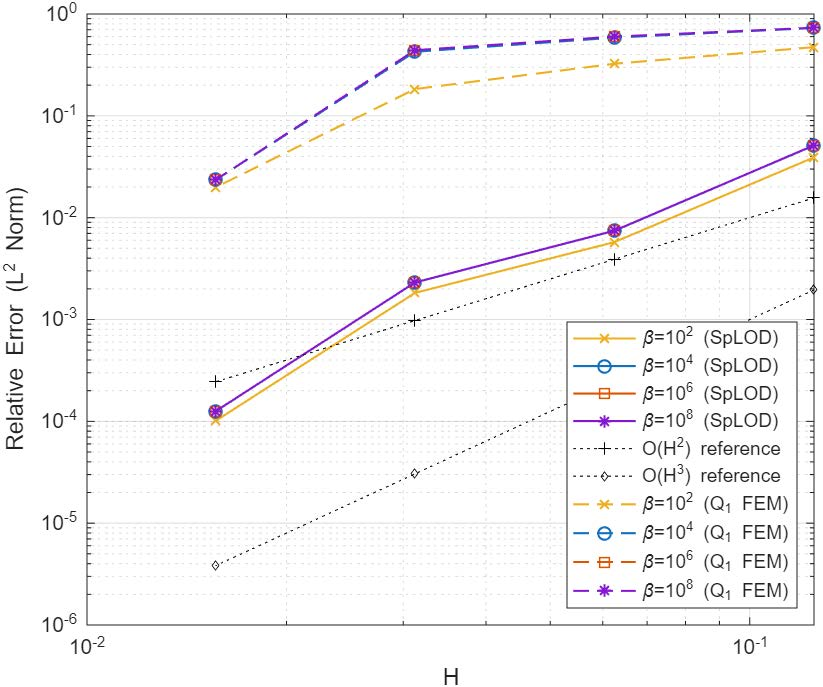}
\end{center}
\caption{Relative energy errors (left) and relative  $L^2$ errors (right) for Example~\ref{example:Irregular}.}
\label{fig:IrregularErrors}
\end{figure}
\par
{The solution times of \eqref{eq:LOD} by the backslash command in MATLAB are
 reported in Table~\ref{tab:IrregularTimes}.}
\begin{table}[H]
\footnotesize
\begin{center}
\begin{tabular}{|c|c|c|c|c|}
\hline
$H \setminus \beta$ &  $10^2$& $10^4$& $10^6$ & $10^8$ \\
\hline
&&&&\\[-11pt]
1/8 &   $8.78\times10^{-5}$ &  $7.47\times10^{-5}$ &  $7.64\times10^{-5}$ &  $7.91\times10^{-5}$\\
\hline
&&&&\\[-11pt]
1/16&   $1.77\times10^{-3}$ &  $1.53\times10^{-3}$&  $1.45\times10^{-3}$&  $1.68\times10^{-3}$\\
\hline
&&&&\\[-11pt]
1/32&   $2.86\times10^{-2}$ &  $2.52\times10^{-2}$ &  $2.28\times10^{-2}$ &  $2.43\times10^{-2}$\\
\hline
\end{tabular}
\end{center}
\caption{{Solution times (in seconds) of \eqref{eq:LOD} by the backslash command} {in MATLAB
for Example~\ref{example:Irregular}}. }
\label{tab:IrregularTimes}
\end{table}
\end{example}
\begin{example}\label{example:GF}
\rm
{The picture in this example depicts a real-life  geographical fracture
 found at
 www.beg.utexas.edu/eichhubl/Pages/Fracturemechanics.html,}
 where the diffusion coefficient
 $\kappa(x)=\beta$ in the black region and $1$ elsewhere
 (cf. Figure~\ref{fig:GF}).
\begin{figure}[htbp]
\begin{center}
  \includegraphics[height=2in]{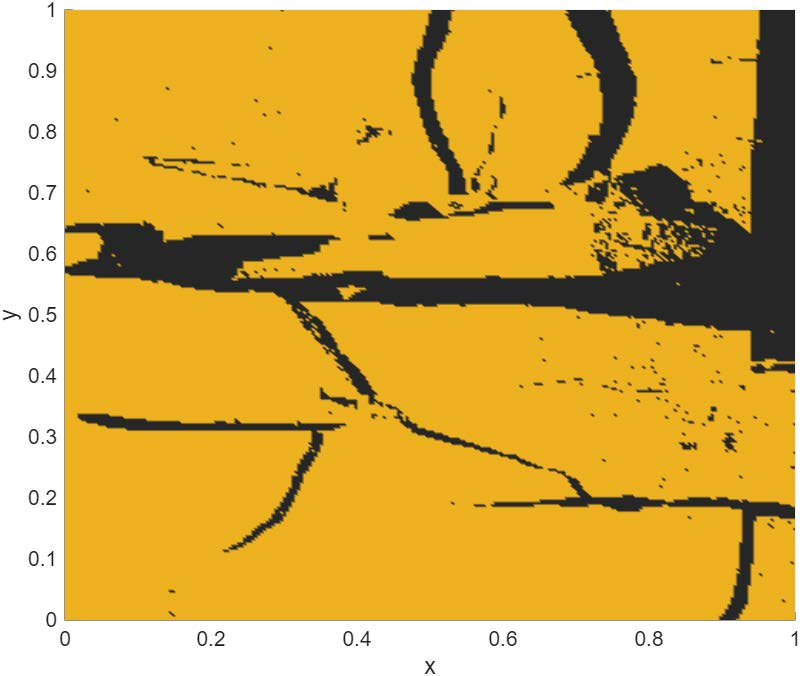}
\end{center}
\par\vspace{-.2in}
\caption{High contrast geographical fracture.}
\label{fig:GF}
\end{figure}
\par
We solve the boundary value problem \eqref{eq:BVP} for $\beta=10^2,10^4,10^6$ and $10^8$.  We take
 $h=1/272$ and $H=1/17, 1/34, 1/68$.
 The data generated by the computation  are displayed in Table~\ref{tab:GFData}.
{
\begin{table}[htbp]
\resizebox{1 \textwidth}{!}{
  \begin{tabular}{|l||c|c|c|c||c|c|c|c||c|c|c|c||c|c|c|c|}
  \hline
  \multirow{2}{*}{\backslashbox{$H$}{$\beta$}}&\multicolumn{4}{c||}{\mystrut(12,4)$10^2$}&\multicolumn{4}{c||}{$10^4$}&
  \multicolumn{4}{c||}{$10^6$}& \multicolumn{4}{c|}{$10^8$}\\
  \cline{2-17}&\mystrut(14,4)$\sqrt{L}$&$\sqrt{M}$&$q$&$k$&$\sqrt{L}$&$\sqrt{M}$&$q$&$k$&
  $\sqrt{L}$&$\sqrt{M}$&$q$&$k$&$\sqrt{L}$&$\sqrt{M}$&$q$&$k$\\
   \hline
   &&&&&&&&&&&&&&&&\\[-9pt]
  $1/17$ &$21$&$3.0\times10^5$&$0.51$&$45$&$22$&$9.7\times10^6$&$0.79$&$129$&$22$&$9.7\times10^8$&$0.80$&$164$
  &$22$&$9.7\times10^{10}$&$0.80$&$195$\\
  \hline&&&&&&&&&&&&&&&&\\[-10pt]
  $1/34$&$38$&$2.9\times10^4$&$0.41$&$27$&$38$&$2.0\times10^5$&$0.41$&$32$&$38$&$8.3\times10^6$&$0.41$&$39$
  &$38$&$5.0\times10^8$&$0.41$&$46$\\
  \hline&&&&&&&&&&&&&&&&\\[-10pt]
  $1/68$&$71$&$1.0\times10^4$&$0.44$&$31$&$71$&$2.0\times10^4$&$0.45$&$35$&$71$&$1.1\times10^5$&$0.45$&$40$
  &$71$&$1.1\times10^6$&$0.45$&$46$\\
  \hline
  \end{tabular}
  }
\par\vspace{.1in}
\caption{Data generated by the computation for Example~\ref{example:GF}.}
\label{tab:GFData}
\end{table}
}
\par
 The estimates for the energy error and $L^2$ error according to
 \eqref{eq:ConcreteEnergyError} and
 \eqref{eq:ConcreteLTError} are shown in Table~\ref{tab:GFEstimates},
 and the actual energy errors and  $L^2$ errors can be found
 in Table~\ref{tab:GFErrors}.
 It is observed that for a given $H$ the errors are independent
 of the contrast $\beta$ and the estimates
 \eqref{eq:ConcreteEnergyError} and
 \eqref{eq:ConcreteLTError} are satisfied.
\begin{table}[htbp]
\begin{center}
{\scriptsize
  \begin{tabular}{|l|c|c|c|}
  \hline&&&\\[-10pt]
    \hspace{52pt}$H$ & $\frac17$ & $\frac{1}{34}$ & $\frac{1}{68}$ \\
    &&&\\[-10pt]
    \hline&&&\\[-10pt]
    Energy Error Estimate&$8.5\times10^{-2}$&$4.3\times10^{-2}$&$2.1\times10^{-2}$\\
    \hline&&&\\[-10pt]
   $L^2$ Error Estimate&$7.3\times10^{-3}$&$1.9\times10^{-3}$&$4.5\times10^{-4}$\\
    \hline
  \end{tabular}
}
\end{center}\par\vspace{.1in}
\caption{Estimates \eqref{eq:ConcreteEnergyError} and
 \eqref{eq:ConcreteLTError} for Example~\ref{example:GF}.}
 \label{tab:GFEstimates}
\end{table}
\begin{table}[htbp]
\resizebox{1 \textwidth}{!}{
  \begin{tabular}{|l||c|c||c|c||c|c||c|c|}
     \hline
  \multirow{2}{*}{\backslashbox{$H$}{$\beta$}}&\multicolumn{2}{c||}{\mystrut(12,4)$10^2$}&\multicolumn{2}{c||}{$10^4$}&
  \multicolumn{2}{c||}{$10^6$}& \multicolumn{2}{c|}{$10^8$}\\
  \cline{2-9}&\mystrut(12,4)Energy Error&$L^2$ Error&Energy Error&$L^2$ Error
  &Energy Error&$L^2$ Error&Energy Error&$L^2$ Error\\
    \hline
   &&&&&&&&\\[-9pt]
  $1/17$&$4.1\times10^{-3}$&$5.1\times10^{-5}$  &$4.2\times10^{-3}$&$5.8\times10^{-5}$
  &$4.2\times10^{-3}$&$5.9\times10^{-5}$   &$4.2\times10^{-3}$& $5.8\times10^{-5}$  \\
  \hline&&&&&&&&\\[-10pt]
  $1/34$&$1.5\times10^{-3}$&$1.0\times10^{-5}$   &$1.4\times10^{-3}$&$1.0\times10^{-5}$
  &$1.4\times10^{-3}$&$1.0\times10^{-5}$
  &$1.4\times10^{-3}$&$1.0\times10^{-5}$  \\
  \hline&&&&&&&&\\[-10pt]
  $1/68$&$4.5\times10^{-4}$&$1.5\times10^{-6}$  &$4.6\times10^{-4}$&$1.6\times10^{-6}$
  &$4.6\times10^{-4}$&$1.6\times10^{-6}$
  &$4.6\times10^{-4}$&$1.6\times10^{-6}$  \\
  \hline
  \end{tabular}
}
\par\vspace{.1in}
\caption{Energy and $L^2$ errors for Example~\ref{example:GF}.}
\label{tab:GFErrors}
\end{table}
\par
 The convergence histories in the energy norm and the  $L^2$ norm for the
 standard $Q_1$ finite element method and the localized
 multiscale finite element method are displayed in Figure~\ref{fig:GFErrors}.
  Their behaviors are similar to those
 in Example~\ref{example:PS} and Example~\ref{example:Irregular}.
\begin{figure}[htbp]
\begin{center}
  \includegraphics[height=2in]{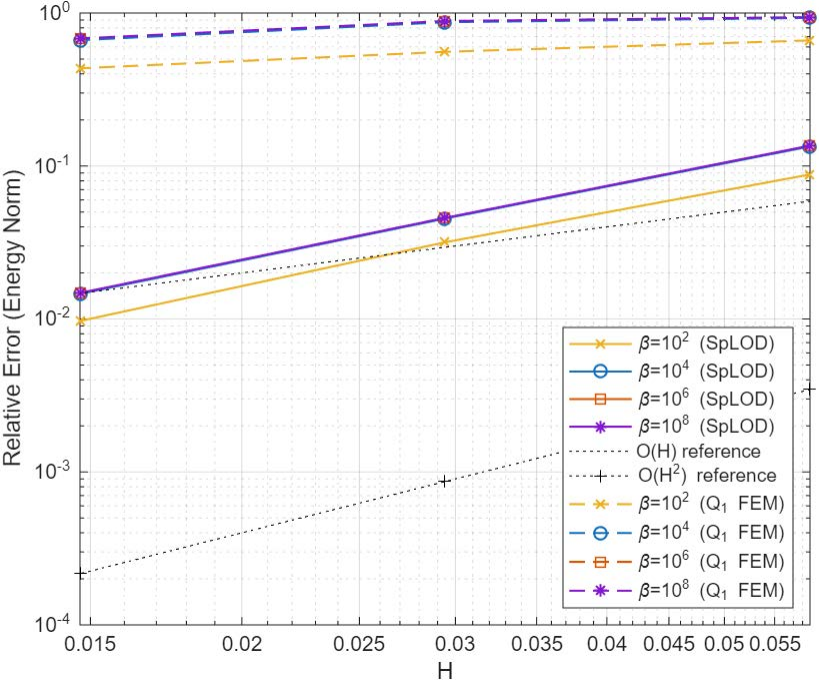} \qquad \qquad
  \includegraphics[height=2in]{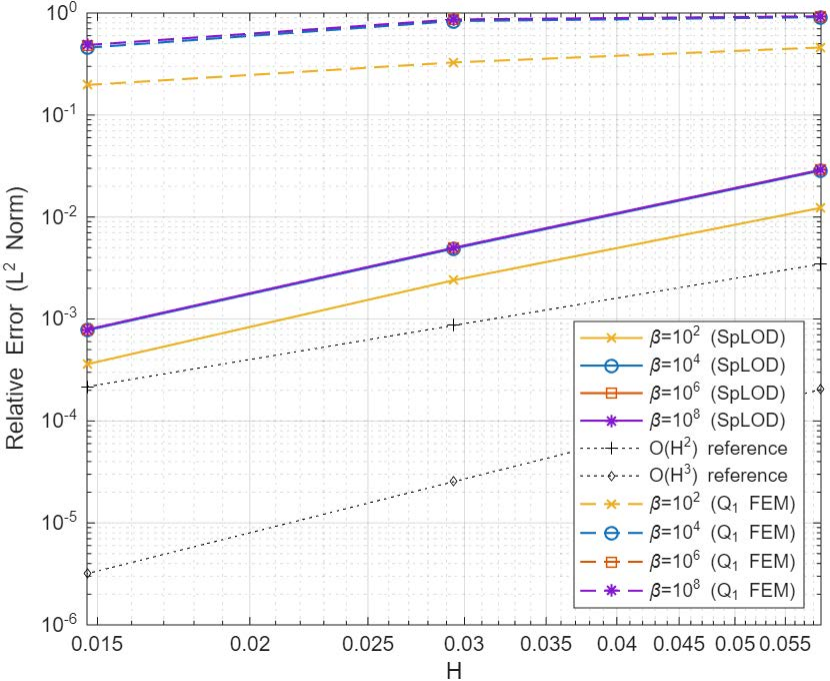}
\end{center}
\caption{Relative energy errors (left) and relative $L^2$ errors (right) for Example~\ref{example:GF}.}
\label{fig:GFErrors}
\end{figure}
\par
{
Table~\ref{tab:GFTimes} contains the solution times of \eqref{eq:LOD}
by the backslash command in MATLAB.}
\begin{table}[H]
\footnotesize
\begin{center}
\begin{tabular}{|c|c|c|c|c|}
\hline
&&&&\\[-11pt]
$H \setminus \beta$ &  $10^2$& $10^4$& $10^6$ & $10^8$ \\
\hline
&&&&\\[-11pt]
1/17 &  $4.95\times10^{-3}$ & $6.48\times10^{-3}$ &  $6.69\times10^{-3}$ &  $7.46\times10^{-3}$\\
\hline
&&&&\\[-11pt]
1/34 &  $4.82\times10^{-2}$  & $5.22\times10^{-2}$ &  $4.71\times10^{-2}$ &  $4.60\times10^{-2}$\\
\hline
\end{tabular}
\end{center}
\caption{{Solution times (in seconds) of \eqref{eq:LOD} by the backslash command}
{in MATLAB for Example~\ref{example:GF}}.}
\label{tab:GFTimes}
\end{table}

\end{example}
\section{Concluding Remarks}\label{sec:Conclusions}
We have developed a multiscale finite element method for a model diffusion problem
with rough and high contrast coefficients whose performance is similar to the performance
of standard finite
element methods for the Poisson equation on smooth or convex domains.
 We established simple explicit error estimates
under conditions that can be verified from the outputs of the computation in the construction of the
multiscale finite element space.
\par
The method was presented for the unit square for clarity.  However,
since the construction and analysis of our method
are based purely on techniques in numerical linear algebra, the extension to arbitrary domains in
two and three dimensions with simplicial triangulations and the $P_1$ finite element
is straightforward.
\par
 It would be interesting to carry out an analysis of the condition number
 of the matrix $\bK^T\bA\bK$ in
 Section~\ref{subsec:Correction} based on domain decomposition techniques
  (cf. Remark~\ref{rem:Preconditioning}).
\appendix
\section{Lower Bounds for the Eigenvalue $\mu_i$}\label{append:EVs}
 We will derive the lower bounds of the eigenvalues stated in Remark~\ref{rem:Scaling}.
\subsection{}\quad
 Let $K_i$ be one of the (closed) elements of $\cT_H$ that is disjoint from $\p\O$.
 By scaling,
 the eigenvalue problem \eqref{eq:EigenEquation}--\eqref{eq:Normalization}
 with $\kappa=1$ is equivalent to
 the following discrete eigenvalue problem on the unit square $\Omega$:
 Find $\mu_j>0$ such that
 \begin{align}
  \int_\Omega \nabla \psi_j\cdot\nabla w\,dx&=
  \mu_j \int_\O
  \psi_j \,w\,dx  \qquad\forall\,w\in V_{h/H},\label{eq:REigenEquation}\\
  \int_\O \psi_j^2dx&=1,\label{eq:RNormalization}
\end{align}
 where $V_{h/H}$ is the $Q_1$ finite element subspace of $H^1(\O)$ associated with a uniform
 triangulation of $\O$ by squares whose edges have length $h/H$.
\par
 Observe that the discrete eigenvalue problem
  \eqref{eq:REigenEquation}--\eqref{eq:RNormalization}
 is a finite element approximation of the
 following continuous eigenvalue problem:  Find $\hat\mu_j>0$ such that
\begin{align}
  \int_\Omega \nabla \hat\psi_j\cdot\nabla  w\,dx&=
  \hat\mu_j \int_\O
 \hat\psi_j \, w\,dx  \qquad\forall\, w\in H^1(\O),\label{eq:RE1}\\
  \int_\O \hat\psi_j^2dx&=1.\label{eq:RN1}
\end{align}
\par
 The smallest eigenvalue for \eqref{eq:RE1}--\eqref{eq:RN1} is given by
$\hat\mu_1=\pi^2$  with eigenfunction $\hat\psi_1(x)=\cos(\pi x_1)$.
It is characterized by
\begin{equation}\label{eq:FirstEigenvalue}
  \hat\mu_1=\min_{\substack{\strut w\in H^1(\O),\, \|w\|_\LT=1\\
  \int_\O w\,dx=0}}\int_\O |\nabla w|^2 dx.
\end{equation}
\par
 On the other hand the smallest eigenvalue $\mu_1$ of
 \eqref{eq:REigenEquation}--\eqref{eq:RNormalization}
 is characterized by
\begin{equation}\label{eq:FirstDE}
  \mu_1=\min_{\substack{\strut w\in V_{h/H},\,
  \|w\|_\LT=1\\ \int_\O w\,dx=0}}\int_\O |\nabla w|^2 dx.
\end{equation}
\par
 Comparing \eqref{eq:FirstEigenvalue} and \eqref{eq:FirstDE}, we see that
\begin{equation}\label{eq:LB1}
  \mu_1\geq\hat\mu_1=\pi^2.
\end{equation}
\subsection{} \quad Let $K_i$ be one of the (closed) elements of
 $\cT_H$ that has only one edge on $\p\O$.
In this case the corresponding continuous eigenvalue problem is to find
$\hat\mu_j$ such that
\begin{align}
  \int_\Omega \nabla \hat\psi_j\cdot\nabla  w\,dx&=
  \hat\mu_j \int_\O
 \hat\psi_j \, w\,dx  \qquad\forall\, w\in H^1(\O),\;w(0,t)=0, \,0<t<1,\label{eq:RE2}\\
  \int_\O \hat\psi_j^2dx&=1.\label{eq:RN2}
\end{align}
\par
 The smallest eigenvalue for \eqref{eq:RE2}--\eqref{eq:RN2} is
 $\hat\mu_1=(\pi/2)^2$ with eigenfunction
 $\hat\psi_1(x)=\sin((\pi/2)x_1)$.  It is characterized by
 $$\hat\mu_1=\min_{\substack{\strut w \in H^1(\O),\,\|w\|_\LT=1\\
   w(0,t)=0, \, 0<t<1}}\int_\O |\nabla w|^2dx$$
 and hence provides a lower bound for the discrete eigenvalue characterized by
  $$\mu_1=\min_{\substack{\strut w \in V_{h/H},\,\|w\|_\LT=1\\
   w(0,t)=0, \, 0<t<1}}\int_\O |\nabla w|^2dx.$$
\subsection{} \quad  Let $K_i$ be one of the (closed) elements of
$\cT_H$ that has two edges on $\p\O$.
In this case the corresponding continuous eigenvalue problem is to find
$\hat\mu_j$ such that
\begin{align}
  \int_\Omega \nabla \hat\psi_j\cdot\nabla  w\,dx&=
  \hat\mu_j \int_\O
 \hat\psi_j \, w\,dx  \qquad\forall\, w\in H^1(\O),\;w(0,t)=w(t,0)=0,\,0<t<1,\label{eq:RE3}\\
  \int_\O \hat\psi_j^2dx&=1.\label{eq:RN3}
\end{align}
\par
\par
 The smallest eigenvalue for \eqref{eq:RE2}--\eqref{eq:RN2} is
 $\hat\mu_1=\pi^2/2$ with eigenfunction
 $$\hat\psi_1(x)=\sin((\pi/2)x_1)\sin((\pi/2)x_2).$$
   It is characterized by
 $$\hat\mu_1=\min_{\substack{\strut w \in H^1(\O),\,\|w\|_\LT=1\\
   w(0,t)=w(t,0)=0, \, 0<t<1}}\int_\O |\nabla w|^2dx$$
 and hence provides a lower bound for the discrete eigenvalue characterized by
  $$\mu_1=\min_{\substack{\strut w \in V_{h/H},\,\|w\|_\LT=1\\
   w(0,t)=w(t,0)=0, \, 0<t<1}}\int_\O |\nabla w|^2dx.$$

\end{document}